\documentclass[reqno]{amsart}

\usepackage[a4paper]{geometry}
\usepackage[headings]{fullpage}
\usepackage[onehalfspacing]{setspace}

\usepackage{graphicx}

\usepackage{mathrsfs}
\usepackage{amsfonts, amssymb, amsmath, amsthm, esint}
\usepackage{cite}
\usepackage{latexsym}
\usepackage{enumerate}
\usepackage{array}
\usepackage{color}
\usepackage{subcaption}
\usepackage{caption}
\usepackage{multirow}

\usepackage{tcolorbox}

\usepackage{cleveref}

\allowdisplaybreaks

\usepackage{tikz}

\allowdisplaybreaks

\theoremstyle{definition}

\newtheorem{assum}{Assumption}

\theoremstyle{theorem}
\newtheorem{theorem}{Theorem}

\newtheorem{lemma}{Lemma}
\newtheorem{spro}{Proposition}

\theoremstyle{remark}
\newtheorem{remark}{Remark}

\newenvironment{eqn*}[0]
{\begin{equation*}}{\end{equation*}\ignorespacesafterend}

\newenvironment{dis}[0]
{\begin{equation*}}{\end{equation*}\ignorespacesafterend}

\newenvironment{eqn}[0]
{\begin{equation}}{\end{equation}\ignorespacesafterend}

\newcommand{\mc}[1]{\mathcal{#1}}

\newcommand{\tr}[1]{\textrm{#1}}

\newcommand{\pd}[0]{\partial}
\newcommand{\ol}[1]{\overline{#1}}

\newcommand{\vphi}[0]{\varphi}

\newcommand{\sus}[0]{\subset}

\newcommand{\valpha}[0]{\boldsymbol{\alpha}}

\newcommand{\vvph}[0]{\boldsymbol{\vphi}}
\newcommand{\vpsi}[0]{\boldsymbol{\psi}}

\newcommand{\vf}[0]{\boldsymbol{f}}
\newcommand{\vg}[0]{\boldsymbol{g}}

\newcommand{\vn}[0]{\boldsymbol{n}}

\newcommand{\vq}[0]{\boldsymbol{q}}

\newcommand{\vt}[0]{\boldsymbol{t}}
\newcommand{\vu}[0]{\boldsymbol{u}}
\newcommand{\vv}[0]{\boldsymbol{v}}
\newcommand{\vw}[0]{\boldsymbol{w}}
\newcommand{\vx}[0]{\boldsymbol{x}}

\newcommand{\zz}[0]{\boldsymbol{0}}

\newcommand{\vH}[0]{\boldsymbol{H}}

\newcommand{\vL}[0]{\boldsymbol{L}}
\newcommand{\vM}[0]{\boldsymbol{M}}

\newcommand{\vP}[0]{\boldsymbol{P}}

\newcommand{\vV}[0]{\boldsymbol{V}}
\newcommand{\vW}[0]{\boldsymbol{W}}

\newcommand{\vZ}[0]{\boldsymbol{Z}}

\newcommand*\diff{\mathop{}\!\mathrm{d}}

\newcommand{\R}[0]{\mathbb{R}}

\newcommand{\abs}[1]{\left|#1\right|}
\newcommand{\enorm}[1]{{\left\vert\kern-0.25ex\left\vert\kern-0.25ex\left\vert #1 
		\right\vert\kern-0.25ex\right\vert\kern-0.25ex\right\vert}}
\newcommand{\inn}[1]{\langle #1\rangle}

\newcommand{\wt}[1]{\widetilde{#1}}

\DeclareMathOperator{\diam}{diam}

\DeclareMathOperator{\rot}{rot}
\DeclareMathOperator{\spn}{span}
\DeclareMathOperator{\curl}{\mathbf{curl}}

\let\div\relax
\DeclareMathOperator{\div}{div}

\title[Divergence-Free Basis in Nonconforming VEM]{A Formal Construction of a Divergence-Free Basis in the Nonconforming Virtual Element Method for the Stokes Problem}

\author{Do Y. Kwak \and Hyeokjoo Park}

\thanks{Department of Mathematical Sciences, Korea Advanced Institute of Science and Technology, Daejeon, 305-701, Korea (kdy@kaist.ac.kr, hjpark235@kaist.ac.kr), This work is partially supported by NRF, contract No. 2021R1A2C1003340.}

\date{\today}

\keywords{Nonconforming virtual element method, Stokes problem, polygonal mesh, divergence-free element}

\subjclass[2010]{65N12, 65N30, 76D07}

\begin{document}

\maketitle

\begin{abstract}
We develop a formal construction of a pointwise divergence-free basis in the nonconforming virtual element method of arbitrary order for the Stokes problem introduced in \cite{zhao2019divergence}. The proposed construction can be seen as a generalization of the divergence-free basis in Crouzeix-Raviart finite element space \cite{thomassetimplementation,brenner1990nonconforming} to the virtual element space. Using the divergence-free basis obtained from our construction, we can eliminate the pressure variable from the mixed system and obtain a symmetric positive definite system. Several numerical tests are presented to confirm the efficiency and the accuracy of our construction.
\end{abstract}

\section{Introduction}

Recently, the virtual element method (VEM) was proposed in \cite{beirao2013basic} as a generalization of the finite element method (FEM) to general polygonal and polyhedral meshes. In VEMs, the local discrete spaces on the mesh polygons/polyhedrons, called local virtual element spaces, consist of polynomials of certain degrees and some other non-polynomial functions that are solutions of specific partial differential equations. Although such functions are not defined explicitly, they are characterized by degrees of freedom, such as values at mesh vertices, the moments on mesh edges/faces, and the moments on mesh polygons/polyhedrons. On each (polygonal or polyhedral) element, the discrete bilinear form can be computed using only the degrees of freedom, and satisfies two properties, called consistency and stability. The consistency means that the discrete bilinear form is equal to the continuous bilinear form when one of the arguments is a polynomial, and the stability means that the discrete bilinear form is coercive for general virtual elements. Moreover, the virtual element spaces can be extended to arbitrary order in straightforward way. Because of such advantages, VEMs have been developed for many different types of equations, and successfully applied to various problems. For more thorough survey, we refer to \cite{beirao2013basic,brezzi2014basic,ahmad2013equivalent,beirao2014hitchhiker,de2016nonconforming,cangiani2017conforming,beirao2016virtual,da2013virtual,zhang2019nonconforming,antonietti2018fully} and references therein. 

There have appeared some results concering the VEMs for the Stokes problem as well. In \cite{antonietti2014stream}, a stream formulation of the VEM for the Stokes problem was presented. In \cite{cangiani2016nonconforming,liu2017nonconforming}, the nonconforming VEM of arbitrary order for the Stokes problem on polygonal and polyhedral meshes was first introduced. Therein, each component of the velocity is approximated by the nonconforming virtual element space presented in \cite{de2016nonconforming}. However, the velocity approximation in \cite{cangiani2016nonconforming,liu2017nonconforming} is not pointwise divergence-free, and it is merely divergence-free in a relaxed (projected) sense. 

In the two-dimensional case, some researchers have developed VEMs for the Stokes problem in which the velocity approximation is pointwise divergence-free. In \cite{da2017divergence}, the divergence-free velocity approximation is presented in the conforming virtual element space of order $k\geq 2$. On each polygon, the virtual element space consists of velocity solutions of the local Stokes problem with Dirichlet boundary condition. On the other hand, the nonconforming virtual element space of arbitrary order was constructed by enriching a $\vH(\div)$-conforming virtual element space in \cite{zhao2019divergence}. However, the proposed methods in \cite{da2017divergence,zhao2019divergence} only showed that the computed velocity approximation is pointwise divergence-free. They do not discuss the construction of divergence-free basis functions. To the best of our knowledge, a formal construction of divergence-free bases in these VEMs has never been considered and developed. 

The main goal of this paper is to present a formal construction of a divergence-free basis in the two-dimensional nonconforming VEMs for the Stokes problem introduced in \cite{zhao2019divergence}. We first compute the dimension of the divergence-free subspace of the nonconforming virtual element space, using Euler's formula. We then construct basis functions of the subspace, in a similar fashion to the divergence-free basis functions proposed in \cite{thomassetimplementation,brenner1990nonconforming} but we generalize to polygonal meshes and higher-order virtual elements. Using the construction of a divergence-free basis, we can eliminate the pressure variable from the coupled system and reduce the saddle point problem to a symmetric positive definite system having fewer unknowns in velocity variable only. Although we only consider the Stokes problem in this paper, we expect that our construction can be applied to more complicated problems, such as the incompressible Navier-Stokes problem.

The rest of this paper is organized as follows. In section 2, we state the stationary Stokes problem and its variational formulation. In section 3, we review the divergence-free nonconforming VEM for the Stokes problem introduced in \cite{zhao2019divergence}. In section 4, we discuss a formal construction of divergence-free basis of the nonconforming virtual element space. In section 5, we discuss implementations including nonhomogeneous Dirichlet boundary conditions. In section 6, we offer some numerical experiments that verify the efficiency and the accuracy of our construction. Finally, conclusions are given in section 7.

\section{Model Problem}

Let $\Omega\sus\R^2$ be a bounded, convex polygonal domain with boundary $\pd\Omega$. We consider the Stokes problem on $\Omega$: Given $\vf:\Omega\to\R^2$ and $\vg:\pd\Omega\to\R^2$, find $\vu:\Omega\to\R^2$ and $p:\Omega\to\R$ such that
\begin{eqn}\label{eqn:Stokes123}
\left\{\begin{array}{rcl}
-\Delta\vu + \nabla p = \vf & \tr{in} & \Omega, \\
\div\vu = 0 & \tr{in} & \Omega, \\
\vu = \vg & \tr{on} & \pd\Omega.
\end{array}\right.
\end{eqn}

In order to obtain the variational formulation of \eqref{eqn:Stokes123}, we introduce the usual notation for Sobolev spaces, norms, seminorms, and inner products. Let $D$ be a bounded domain in $\R^2$. We then define $\vL^2(D) = [L^2(D)]^2$ and $\vH^s(D) = [H^s(D)]^2$ for $s > 0$. The $L^2$-inner product of $L^2(D)$ and $\vL^2(D)$ is denoted by $(\cdot,\cdot)_{0,D}$. Next, for $s \geq 0$, the $H^s$-norm of $H^s(D)$ and $\vH^s(D)$ is denoted by $\|\cdot\|_{s,D}$. Similarly, for $s > 0$, the $H^s$-seminorm of $H^s(D)$ and $\vH^s(D)$ is denoted by $|\cdot|_{s,D}$. The subspace $L_0^2(D)$ of $L^2(D)$ is defined by
\begin{dis}
L_0^2(D) = \left\{q\in L^2(D) : \int_Dq\diff\vx = 0\right\}.
\end{dis}
Let us define
\begin{eqnarray*}
\vH_0^1(\Omega) & = & \left\{\vv\in \vH^1(\Omega) : \vv = \zz \ \tr{on} \ \pd\Omega\right\}, \\
\vH_{\vg}^1(\Omega) & = & \left\{\vv\in\vH^1(\Omega) : \vv = \vg \ \tr{on} \ \pd\Omega\right\}.
\end{eqnarray*}
Then the variational form of the Stokes problem \eqref{eqn:Stokes123} is written as follows: For a given $\vf\in \vL^2(\Omega)$ and a given $\vg\in \vH^{1/2}(\pd\Omega)$ satisfying
\begin{eqn}\label{eqn:StokesBdry}
\int_{\pd\Omega}\vg\cdot\vn_{\Omega}\diff s = 0
\end{eqn}
where $\vn_{\Omega}$ is the unit normal vector on $\pd\Omega$ in the outward direction with respect to $\Omega$, find $\vu\in\vH_{\vg}^1(\Omega)$ and $p\in L_0^2(\Omega)$ such that
\begin{eqn}\label{eqn:StokesVarCon}
\left\{\begin{array}{rcll}
a(\vu,\vv) + b(\vv,p) & = & (\vf,\vv)_{0,\Omega} & \forall \vv\in \vH_0^1(\Omega), \\
b(\vu,q) & = & 0 & \forall q\in L_0^2(\Omega),
\end{array}\right.
\end{eqn}
where 
\begin{eqn}\label{eqn:StokesVarConBilinear}
a(\vu,\vv) = \int_{\Omega}\nabla\vv:\nabla\vu\diff\vx, \quad b(\vv,q) = -\int_{\Omega}q\div\vv\diff\vx.
\end{eqn}
The functions $\vu$ and $p$ are called velocity and pressure, respectively. 


\section{Divergence-Free Nonconforming VEM for the Stokes Problem}\label{sec:DivFreeNCVEM}

In this section, we summarize some preliminaries and review the divergence-free nonconforming VEM for the Stokes problem introduced in \cite{zhao2019divergence}.

\subsection{Notations and preliminaries}

Let $\{\mc{P}_h\}_h$ be a family of decompositions (meshes) of the domain $\Omega$ into polygonal elements $K$ with maximum diameter $h$. We assume that the decompositions satisfy the following regularity properties \cite{beirao2013basic,de2016nonconforming,da2017divergence,zhao2019divergence}.
\begin{assum}
There exists $\rho > 0$ independent of $h$ such that
\begin{itemize}
\item the decomposition $\mc{P}_h$ consists of a finite number of nonoverlapping convex polygonal elements;
\item if $K\in\mc{P}_h$ and $e$ is an edge of $K$ then $h_e \geq \rho h_K$, where $h_e$ and $h_K$ denote the diameter of $e$ and $K$, respectively;
\item every element $K$ of $\mc{P}_h$ is star-shaped with respect to the ball of radius $\rho h_K$.
\end{itemize} 
\end{assum}
We next define some notations for sets of mesh items. We denote by $\mc{V}_{h}$ and $\mc{E}_h$ the set of all mesh vertices and mesh edges in $\mc{P}_h$, respectively. We also denote by $\mc{V}^i_h$ and $\mc{V}^\pd_h$ the set of all mesh vertices in the internal and the boundary of $\mc{P}_h$, respectively. Similarly $\mc{E}_h^i$ is the set of all mesh edges in the internal of $\mc{P}_h$, and $\mc{E}_h^\pd$ the set of all mesh edges in the boundary of $\mc{P}_h$. We also define
\begin{eqnarray*}
N_P & = & \tr{the number of polygons in $\mc{P}_h$}, \\
N_E & = & \tr{the number of edges in $\mc{E}_h$}, \\
N_{E,i} & = & \tr{the number of edges in $\mc{E}_h^i$}, \\
N_{E,\pd} & = & \tr{the number of edges in $\mc{E}_h^{\pd}$}, \\
N_V & = & \tr{the number of vertices in $\mc{V}_h$}, \\
N_{V,i} & = & \tr{the number of vertices in $\mc{V}_h^i$}, \\
N_{V,\pd} & = & \tr{the number of vertices in $\mc{V}_h^{\pd}$}.
\end{eqnarray*}
For each $K\in\mc{P}_h$, let $\vn_K$ and $\vt_K$ denote its exterior unit normal vector and counterclockwise tangential vector, respectively. Let $e\in\mc{E}_h^i$. We then define respectively $\vn_e$ and $\vt_e$ as the unit normal and tangential vector of $e$ with orientation fixed once and for all. Next let $e\in\mc{E}_h^{\pd}$, we define respectively $\vn_e$ and $\vt_e$ as the unit normal and tangential vector on $e$ in the outward and counterclockwise direction with respect to $\Omega$.

Let $e\in\mc{E}_h^i$ and let $K^-$ and $K^+$ be the polygons in $\mc{P}_h$ that have $e$ as a common edge, and satisfy $\vn_e = \vn_{K_+}$ on $e$ (i.e., $\vn_e$ points from $K^+$ to $K^-$). If $e\in\mc{E}_h^{\pd}$, we define $\vn_e$ by the unit normal vector in the outward direction with respect to $\Omega$.

Again let $e\in\mc{E}_h^i$ and let $K^-$ and $K^+$ are the polygons in $\mc{P}_h$ having $e$ as a common edge. For $\vv:\Omega\to\R^2$ satisfying $\vv|_{K^+}\in \vH^1(K^+)$ and $\vv|_{K^-}\in \vH^1(K^-)$, we define the jump of $\vv$ on the edge $e$ by
\begin{dis}
[\vv]_e = \vv|_{K^+}(\vn_e\cdot\vn_{K^+}) + \vv|_{K^-}(\vn_e\cdot\vn_{K^-}).
\end{dis}
If $e\in\mc{E}_h^{\pd}$, we define $[\vv]_e = \vv|_e$. 

We define the broken Sobolev space $\vH^1(\Omega;\mc{P}_h)$ by
\begin{dis}
\vH^1(\Omega;\mc{P}_h) = \left\{\vv\in\vL^2(\Omega) : \vv|_K\in \vH^1(K) \ \forall K\in\mc{P}_h\right\}
\end{dis}
and define its norm and seminorm by
\begin{dis}
\|\vv\|_{1,h} = \left(\sum_{K\in\mc{P}_h}\|\vv\|_{1,K}^2\right)^{1/2}, \quad |\vv|_{1,h} = \left(\sum_{K\in\mc{P}_h}|\vv|_{1,K}^2\right)^{1/2}.
\end{dis}
We also define
\begin{dis}
\vH^{1,nc}(\Omega;\mc{P}_h) = \left\{\vv\in\vH^1(\Omega;\mc{P}_h) : \int_e[\vv]_e\cdot\vq\diff\vx = 0 \ \forall \vq\in \vP_{k-1}(e), \ \forall e\in\mc{E}_h^i\right\}.
\end{dis}

Let $O$ be an $1$ or $2$ dimensional geometrical object (edge or polygon). For an integer $k\geq 0$, $P_k(O)$ denotes the space of polynomials of degree $\leq k$ on $O$. $M_k(O)$ denotes the set of scaled monomials of degree $\leq k$ on $O$, that is,
\begin{dis}
M_k(O) = \left\{\left(\frac{\vx - \vx_O}{h_O}\right)^{\valpha}: |\valpha| \leq k \right\},
\end{dis}
where $\vx$ is a local coordinate system on $O$, $\vx_O$ is the barycenter of $O$ in the local coordinate system, $\valpha$ is a multi-index, and $h_O = \diam(O)$. 

Conventionally we define $P_{-1}(O) = \{0\}$. We also define $\vP_k(O) = (P_k(O))^2$ for $k\geq -1$ and $\vM_k(O) = (M_k(O))^2$ for any nonnegative integer $k$. 

Let $K\in\mc{P}_h$ and let $k$ be a nonnegative integer. We define $(\nabla P_{k+1}(K))^\oplus$ as the subspace of $\vP_k(K)$ satisfying 
\begin{dis}
\vP_k(K) = \nabla P_{k+1}(K) \oplus (\nabla P_{k+1}(K))^\oplus,
\end{dis}
and denote by $\vM_k^\oplus$ a basis of the space $(\nabla P_{k+1}(K))^\oplus$.
For example, one can choose
\begin{dis}
(\nabla P_{k+1}(K))^\oplus = \vx^\perp P_{k-1}(K), \quad \vM_k^\oplus = \left\{m(\vx)\vx^{\perp} : m\in M_{k-1}(K)\right\},
\end{dis}
where $\vx^{\perp} = (x_2,-x_1)$ with $\vx = (x_1,x_2)$.

\subsection{Virtual element space}

We first define a local virtual element space on each element $K\in\mc{P}_h$. Let $k$ be a fixed positive integer. Let 
\begin{dis}
\vW_h^1(K) := \left\{\vv\in \vH^1(K) : \div\vv\in P_{k-1}(K), \ \rot\vv = 0, \ \vv\cdot\vn_K|_e\in P_k(e), \ \forall e\sus\pd K\right\},
\end{dis}
where $\rot\vv = \frac{\pd v_1}{\pd x_2} - \frac{\pd v_2}{\pd x_1}$ for $\vv = (v_1,v_2)\in\vH^1(K)$. Also, let
\begin{dis}
\Phi_h(K) := \left\{\phi\in H^2(K) : \Delta^2\phi\in P_{k-3}(K), \ \phi|_e = 0, \ \Delta\phi|_e\in P_{k-1}(e), \ \forall e\sus\pd K\right\}
\end{dis}
with the convention that $P_{-1}(K) = P_{-2}(K) = \{0\}$. In \cite[Lemma 2]{zhao2019divergence}, it was shown that $\vW_h^1(K) \cap \curl\Phi_h(K) = \{0\}$, where $\curl q = (-\frac{\pd q}{\pd x_2},\frac{\pd q}{\pd x_1})$ for $q\in H^1(K)$. It was also shown in \cite[Lemma 3]{zhao2019divergence} that if the local space $\tilde{\vV}_h(K)$ is defined by
\begin{dis}
\tilde{\vV}_h(K) = \vW_h^1(K) \oplus \curl\Phi_h(K),
\end{dis}
then the following degrees of freedom (DOFs) are unisolvent for $\tilde{\vV}_h(K)$:
\begin{eqnarray*}
\tr{the moments} \ \frac{1}{|e|}\int_e\vv\cdot\vn_eq\diff s, & \quad & q \in M_{k}(e), \\
\tr{the moments} \ \frac{1}{|e|}\int_e\vv\cdot\vt_eq\diff s, & \quad & q \in M_{k-1}(e), \\
\tr{the moments} \ \frac{1}{|K|}\int_K\vv\cdot\vq\diff\vx, & \quad & \vq \in \vM_{k-2}(K).
\end{eqnarray*}
We define a local projection $\Pi^{\nabla}_K:\vH^1(K)\to \vP_k(K)$ on each polygon $K$ in $\mc{P}_h$. It is defined by
\begin{eqnarray*}
\int_K\nabla(\Pi^{\nabla}_K\vv):\nabla\vq\diff\vx & = & \int_K\nabla\vv:\nabla\vq\diff\vx, \quad \forall \vv\in\vH^1(K), \ \forall \vq\in\vP_k(K), \\
\int_{\pd K}\Pi^{\nabla}_K\vv\diff s & = & \int_{\pd K}\vv\diff s,
\end{eqnarray*}
for $\vv\in \vH^1(K)$. Note that $\Pi_K^\nabla\vq = \vq$ for any $\vq\in \vP_k(K)$ and the local projection $\Pi_K^\nabla$ is computable using only the moments of $\vv$ up to order $(k-1)$ on each edge $e\sus \pd K$ and the moments of $\vv$ up to order $(k-2)$ on $K$. 

Now the local nonconforming virtual element space $\vV_h(K)$ on $K$ is defined by
\begin{dis}
\vV_h(K) = \left\{\vv\in\tilde{\vV}_h(K) : \int_e(\vv - \Pi^{\nabla}_K\vv)\cdot\vn_eq\diff s = 0, \ \forall q\in P_k(e)/P_{k-1}(e), \ \forall e\sus \pd K\right\},
\end{dis}
where $P_k(e)/P_{k-1}(e)$ is the subspace of polynomials in $P_k(e)$ that are $L^2(e)$-orthogonal to $P_{k-1}(e)$. It was shown in \cite{zhao2019divergence} that the following DOFs are unisolvent for $\vV_h(K)$:
\begin{eqnarray*}
\frac{1}{|e|}\int_e\vv\cdot\vn_eq\diff s, & \quad & q \in M_{k-1}(e), \\
\frac{1}{|e|}\int_e\vv\cdot\vt_eq\diff s, & \quad & q \in M_{k-1}(e), \\
\frac{1}{|K|}\int_K\vv\cdot\vq\diff\vx, & \quad & \vq \in \vM_{k-2}(K).
\end{eqnarray*}
For each $i = 1,2,\cdots,N_K := \dim\vV_h(K)$, let $\chi_i$ be the operator associated to the $i$-th local DOF. Then for any $\vv\in\vH^1(K)$ there exists a unique element $I_h^K\vv\in\vV_h(K)$ such that 
\begin{dis}
\chi_i(\vv - I_h^K\vv) = 0 \quad \forall i = 1,2,\cdots,N_K.
\end{dis}
The operator $\vv\mapsto I_h^K\vv$ is called a local interpolation operator for $\vV_h(K)$. It was shown in \cite{zhao2019divergence} that we can obtain the following interpolation error estimates.

\begin{spro}[see {\cite[Lemma 6]{zhao2019divergence}}]
There exists a positive constant $C$ independent of $h$ such that for every $K\in\mc{P}_h$ and every $\vv\in\vH^s(K)$ with $1\leq s\leq k+1$,
\begin{dis}
\|\vv - I_h^K\vv\|_{0,K} + h|\vv - I_h^K\vv|_{1,K} \leq Ch^{s}|\vv|_{s,K}.
\end{dis}
\end{spro}

The global nonconforming virtual element spaces are defined as follows:
\begin{eqnarray*}
\vV_h & = & \left\{\vv_h\in \vL^2(\Omega) : \vv_h|_K\in \vV_h(K) \quad \forall K\in\mc{P}_h, \ \int_e[\vv_h]_e\cdot\vq\diff s = 0 \quad \forall \vq\in \vP_{k-1}(e), \ \forall e\in\mc{E}_h^i\right\}, \\
\vV_{h,0} & = & \left\{\vv_h\in \vL^2(\Omega) : \vv_h|_K\in \vV_h(K) \quad \forall K\in\mc{P}_h, \ \int_e[\vv_h]_e\cdot\vq\diff s = 0 \quad \forall \vq\in \vP_{k-1}(e), \ \forall e\in\mc{E}_h\right\}.
\end{eqnarray*}
The global DOFs for $\vV_h$ can be chosen as, for any edge $e$ and polygon $K$ in $\mc{P}_h$,
\begin{eqnarray}
\chi_{e,q}^n(\vv_h) := \frac{1}{|e|}\int_e\vv_h\cdot\vn_eq\diff s, & \quad & q \in M_{k-1}(e), \label{eqn:EdgeNDOFG} \\
\chi_{e,q}^t(\vv_h) := \frac{1}{|e|}\int_e\vv_h\cdot\vt_eq\diff s, & \quad & q \in M_{k-1}(e), \label{eqn:EdgeTDOFG} \\
\chi_{K,\vq}(\vv_h) := \frac{1}{|K|}\int_K\vv_h\cdot\vq\diff\vx, & \quad & \vq \in \nabla M_{k-1}(K) + \vM_{k-2}^{\oplus}(K). \label{eqn:CellDOFG}
\end{eqnarray}
Similarly, the global DOFs for $\vV_{h,0}$ can be chosen. We also define the global interpolation operator $I_h:\vH^{1,nc}(\Omega;\mc{P}_h)\to \vV_h$ by $(I_h\vv)|_K = I_h^K(\vv|_K)$ for each $K\in\mc{P}_h$ and $\vv\in\vH^{1,nc}(\Omega;\mc{P}_h)$.

The discrete pressure space $Q_h$ is defined by
\begin{dis}
Q_h = \{q_h\in L^2_0(\Omega) : q_h|_K\in P_{k-1}(K) \ \forall K\in\mc{P}_h\}.
\end{dis}
The global DOFs for the space $Q_h$ can be chosen as
\begin{dis}
\frac{1}{|K|}\int_Kq_h\phi \diff\vx, \quad \phi\in M_{k-1}(K), \ K\in\mc{P}_h.
\end{dis}
It was shown in \cite{zhao2019divergence} that $\div\vV_h(K)\sus P_{k-1}(K)$ for each $K\in\mc{P}_h$, and $\div_h\vV_{h,0}\sus Q_h$, where $\div_h$ denotes the discrete divergence operator defined by $(\div_h\vv_h)|_K = \div(\vv_h|_K)$ for each $K\in\mc{P}_h$ and $\vv_h\in\vV_h$. Therefore, the nonconforming virtual element space $\vV_h$ is divergence-free. 

\subsection{The discrete problem}

We define a local discrete bilinear form $a_h^K$ for each polygon $K$ in $\mc{P}_h$, as follows.
\begin{dis}
a_h^K(\vv_h,\vw_h) = a^K(\Pi_K^\nabla(\vv_h),\Pi_K^\nabla(\vw_h)) + S^K((I - \Pi_K^\nabla)\vv_h,(I - \Pi_K^\nabla)\vw_h), \quad \vv_h,\vw_h\in \vV_h(K),
\end{dis}
where $a^K$ is the bilinear form defined by
\begin{dis}
a^K(\vv,\vw) = \int_K\nabla\vv:\nabla\vw\diff\vx, \quad \vv,\vw\in\vH^1(K),
\end{dis}
and $S^K$ is a symmetric positive definite bilinear form defined as
\begin{dis}
S^K(\vv_h,\vw_h) = \sum_{i=1}^{N_K}\chi_i(\vv_h)\chi_i(\vw_h), \quad \vv_h,\vw_h\in \vV_h(K),
\end{dis}
where $N_K = \dim(\vV_h(K))$ and $\chi_i$ denotes the operator associated to the $i$-th local DOF for $i = 1,2,\cdots,N_K$. As described in \cite{beirao2013basic,zhao2019divergence}, we obtain the $k$-consistency and stability of $a_h^K$:
\begin{itemize}
\item ($k$-consistency) $a_h^K(\vq,\vv_h) = a^K(\vq,\vv_h)$ for any $\vq\in \vP_k(K)$, $\vv_h\in\vV_h(K)$;
\item (stability) there exist constants $c_*,c^*>0$ independent of $h$ such that
\begin{dis}
c_*a^K(\vv_h,\vv_h) \leq a_h^K(\vv_h,\vv_h) \leq c^*a^K(\vv_h,\vv_h) \quad \forall \vv_h\in \vV_h(K).
\end{dis}
\end{itemize}

The global bilinear form $a_h$ is defined by
\begin{dis}
a_h(\vv_h,\vw_h) = \sum_{K\in\mc{P}_h}a_h^K(\vv_h,\vw_h), \quad \vv_h,\vw_h\in \vV_h.
\end{dis}
On the other hand, the discrete bilinear form $b_h$ is simply defined by
\begin{dis}
b_h(\vv_h,q_h) = \sum_{K\in\mc{P}_h}b^K(\vv_h,q_h), \quad \vv_h\in\vV_h,\ q_h\in Q_h,
\end{dis}
where
\begin{dis}
b^K(\vv,q) = -\int_Kq\div\vv\diff\vx,
\end{dis}
for $\vv\in\vH^1(K)$, $q\in P_{k-1}(K)$, and $K\in\mc{P}_h$. Note that $b_h(\vv_h,q_h)$ is also computable using only the DOFs \eqref{eqn:EdgeNDOFG}-\eqref{eqn:CellDOFG} and we do not rely on the discrete version of it, indeed we omit the subscript $h$ on such bilinear form.

We next discretize the right-hand side $(\vf,\cdot)_{0,\Omega}$ as follows:
\begin{dis}
\inn{\vf_h,\vv_h} = \left\{\begin{array}{ll}
(\vf_h,\ol{\vv}_h)_{0,\Omega} & \tr{if $k = 1$} \\
(\vf_h,\vv_h)_{0,\Omega} & \tr{if $k > 1$}
\end{array}\right., \quad \vv_h\in\vV_h,
\end{dis}
where $\vf_h,\ol{\vv}_h\in \vL^2(\Omega)$ are defined by
\begin{dis}
\vf_h|_K = \left\{\begin{array}{ll}
\Pi_0^K\vf & \tr{if $k = 1$} \\
\Pi_{k-2}^K\vf & \tr{if $k > 1$}
\end{array}\right., \quad \ol{\vv}_h|_K = \frac{1}{|\pd K|}\int_{\pd K}\vv_h\diff s, \quad K\in\mc{P}_h.
\end{dis}
Here, $\Pi_{\ell}^K$ denotes the $L^2$-projection operator onto $\vP_{\ell}(K)$ for each $K\in\mc{P}_h$.

In order to consider the nonhomogeneous Dirichlet boundary condition, let
\begin{dis}
\vV_{h,\vg} = \left\{\vv_h\in\vV_h : \int_e\vg\cdot\vq\diff s = \int_e\vv_h\cdot\vq\diff s, \ \forall \vq\in \vP_{k-1}(e), \ \forall e\in\mc{E}_h^{\pd}\right\}.
\end{dis}
We formulate the nonconforming VEM for the Stokes problem \eqref{eqn:StokesVarCon} as follows: Find $\vu_h\in \vV_{h,\vg}$ and $p_h\in Q_h$ such that
\begin{eqn}\label{eqn:StokesDis}
\left\{\begin{array}{rcll}
a_h(\vu_h,\vv_h) + b_h(\vv_h,p_h) & = & \inn{\vf_h,\vv_h} & \forall \vv_h\in\vV_{h,0}, \\
b_h(\vu_h,q_h) & = & 0 & \forall q_h\in Q_h.
\end{array}\right.
\end{eqn}
Here $\vu_h$ and $p_h$ will be called discrete velocity and discrete pressure, respectively. It was shown in \cite{zhao2019divergence} that the discrete problem \eqref{eqn:StokesDis} is well-posed. Moreover, for the case $\vg = \zz$, we can obtain the following error estimate.

\begin{theorem}[see {\cite[Theorem 13]{zhao2019divergence}}]\label{thm:Error}
Suppose that $\vf\in \vH^{k-1}(\Omega)$ and $\vg = \zz$. Let $(\vu,p)\in (\vH_{0}^1(\Omega)\cap \vH^{k+1}(\Omega))\times(L_0^2(\Omega)\cap H^k(\Omega))$ be the solution of the continuous problem \eqref{eqn:StokesVarCon}. Let $(\vu_h,p_h)\in\vV_{h,0}\times Q_h$ be the solution of the discrete problem \eqref{eqn:StokesDis}. Then 
\begin{dis}
\abs{\vu - \vu_h}_{1,h} + \|p - p_h\|_{0,\Omega} \leq Ch^k\left(|\vu|_{k+1,\Omega} + |p|_{k,\Omega} + |\vf|_{k-1,\Omega}\right),
\end{dis}
where $C$ is a positive constant independent on $h$.  
\end{theorem}

\section{A Formal Construction of Divergence-Free Basis}\label{sec:NCVEM1}

In this section, we present a formal construction of a divergence-free basis for the virtual element space $\vV_{h,0}$. 

We first define the canonical basis associated with the DOFs \eqref{eqn:EdgeNDOFG}-\eqref{eqn:CellDOFG} of the space $\vV_h$. Recall that the global DOFs of $\vV_h$ are given by
\begin{eqnarray*}
\chi_{e,q}^n(\vv_h) = \frac{1}{|e|}\int_e\vv_h\cdot\vn_eq\diff s, & \quad & q \in M_{k-1}(e), \ e\in\mc{E}_h, \\
\chi_{e,q}^t(\vv_h) = \frac{1}{|e|}\int_e\vv_h\cdot\vt_eq\diff s, & \quad & q \in M_{k-1}(e), \ e\in\mc{E}_h, \\
\chi_{K,\vq}(\vv_h) = \frac{1}{|K|}\int_K\vv_h\cdot\vq\diff\vx, & \quad & \vq \in \nabla M_{k-1}(K) + \vM_{k-2}^{\oplus}(K), \ K\in\mc{P}_h.
\end{eqnarray*}
We sometimes write $\chi$ to denote $\chi_{e,q}^n$, $\chi_{e,q}^t$, or $\chi_{K,\vq}$ when it is clear from the context. Using these notations, we define the canonical basis functions of $\vV_h$ associated to the DOFs \eqref{eqn:EdgeNDOFG}-\eqref{eqn:CellDOFG} as follows:
\begin{itemize}
\item For $e\in\mc{E}_h$ and $q\in M_{k-1}(e)$, let $\vvph_{e,q}^n$ be the function in $\vV_h$ such that $\chi_{e,q}^n(\vvph_{e,q}^n) = 1$ and $\chi(\vvph_{e,q}^n) = 0$ for all other DOFs.
\item For $e\in\mc{E}_h$ and $q\in M_{k-1}(e)$, let $\vvph_{e,q}^t$ be the function in $\vV_h$ such that $\chi_{e,q}^t(\vvph_{e,q}^t) = 1$ and $\chi(\vvph_{e,q}^t) = 0$ for all other DOFs.
\item For $K\in\mc{P}_h$ and $\vq\in(\nabla M_{k-1}(K)) + \vM_{k-2}^{\oplus}(K)$, let $\vvph_{K,\vq}$ be the function in $\vV_h$ such that $\chi_{K,\vq}(\vvph_{K,\vq}) = 1$ and $\chi = 0$ for all other DOFs. 
\end{itemize}

Let us define
\begin{dis}
\vZ_{h} = \left\{\vv_h\in\vV_{h} : \div_h\vv_h = 0\right\}, \quad \vZ_{h,0} = \left\{\vv_h\in\vV_{h,0} : \div_h\vv_h = 0\right\}.
\end{dis}
We first compute the dimension of $\vZ_{h,0}$. 

\begin{spro}\label{prop:DimZh}
The dimension of $\vZ_{h,0}$ is 
\begin{dis}
N_{V,i} + kN_{E,i} + (k-1)N_{E,i} + \frac{(k-1)(k-2)}{2}N_{P}.
\end{dis}
\end{spro}

\begin{proof}
Since $\div_h\vV_{h,0} = Q_h$ and since $\div_h\vV_{h,0} \cong \vV_{h,0}/\vZ_{h,0}$, we obtain
\begin{dis}
\dim\vZ_{h,0} = \dim\vV_{h,0} - \dim\left(\div_h\vV_{h,0}\right) = \dim\vV_{h,0} - \dim Q_h.
\end{dis}
Note that
\begin{dis}
\dim\vV_{h,0} = 2\left(\frac{k(k-1)}{2}N_{P} + kN_{E,i}\right), \quad \dim Q_{h} = \frac{k(k+1)}{2}N_P - 1,
\end{dis}
Since $N_{P} - N_{E,i} + N_{V,i} = 1$ from Euler's formula, we have
\begin{eqnarray*}
\dim\vZ_{h,0} & = & \dim\vV_{h,0} - \dim Q_{h} = k(k-1)N_{P} + 2kN_{E,i} - \frac{k(k+1)}{2}N_{P} + 1 \\
& = & k(k-1)N_{P} + 2kN_{E,i} - \frac{k(k+1)}{2}N_{P} + N_{P} - N_{E,i} + N_{V,i} \\
& = & N_{V,i} + kN_{E,i} + (k-1)N_{E,i} + \frac{(k-1)(k-2)}{2}N_{P}.
\end{eqnarray*}
This concludes the proof of the proposition.
\end{proof}

\begin{figure}[!ht]
\centering
\begin{minipage}[b]{.4\linewidth}
\begin{center}
\begin{tikzpicture}
\draw[thick] (0,0) -- (0.4,-0.8) -- (1,-1) -- (2,-0.5) -- (1.5,0.5) -- (0.9,0.8) -- (0,0);
\draw[thick] (0.9,0.8) -- (0.5,1.4) -- (-0.5,1.4) -- (-1,0.8) -- (-1.2,0.2) -- (0,0);
\draw[thick] (-1.2,0.2) -- (-1.5,-0.1) -- (-1.5,-0.8) -- (-1.1,-1.2) -- (-0.6,-1.2) -- (0.4,-0.8);

\draw[thick, ->](0.2,-0.4) -- (0.6,-0.2);
\draw[thick, ->](-0.6,0.1) -- (-0.675,-0.35);
\draw[thick, ->](0.45,0.4) -- (0.05,0.85);

\fill (0,0) circle [radius = 0.075];

\node at (-0.1,0.2) {$v$};

\draw (0.966,-0.166) circle (1.5pt);
\draw (-0.15,0.633) circle (1.5pt);
\draw (-0.6,-0.55) circle (1.5pt);

\node at (0.1,-0.6){\footnotesize{$e_1$}};
\node at (0.6,0.25){\footnotesize{$e_2$}};
\node at (-0.5,0.25){\footnotesize{$e_3$}};

\node at (0.75,-0.55){\footnotesize{$\vn_{e_1,v}$}};
\node at (0.05,1){\footnotesize{$\vn_{e_2,v}$}};
\node at (-1.07,-0.3){\footnotesize{$\vn_{e_3,v}$}};
\end{tikzpicture}
\end{center}
\subcaption{}
\label{fig:ConstructD1}
\end{minipage}
\begin{minipage}[b]{.4\linewidth}
\begin{center}
\begin{tikzpicture}
\draw[thick] (0,0) -- (0.4,-0.8) -- (1.5,-1) -- (2,0.3) -- (0.9,0.8) -- (0,0);
\draw[thick] (0,0) -- (-1,0.2) -- (-1.4,-0.6) -- (-1.2,-1.2) -- (-0.5,-1.5) -- (0.4,-0.8) -- (0,0);

\draw[thick, ->] (0.2,-0.4) -- (0.6,-0.2);

\node at (0.0,-0.5) {$e$};
\node at (0.7,-0.5) {\footnotesize{$\vn_{e}$}};

\draw (0.96,-0.14) circle (1.5pt);
\draw (-0.62,-0.65) circle (1.5pt);
\end{tikzpicture}
\end{center}
\subcaption{}
\label{fig:ConstructD3}
\end{minipage}
\caption{Examples of functions defined in (D1) (left) and (D3) (right).}
\label{fig:ConstructD1D3}
\end{figure}
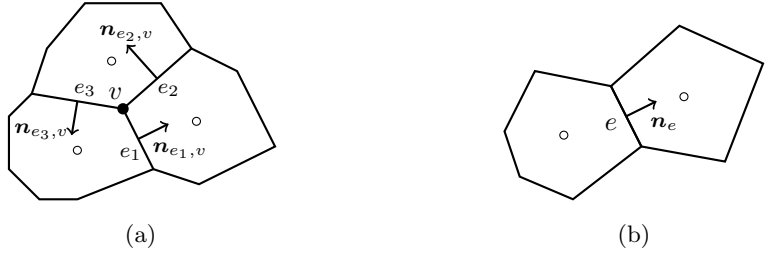

In order to construct a basis of $\vZ_{h,0}$, we first define some functions in $\vV_{h}$. 
\begin{enumerate}[\bfseries\sffamily (D1)]
\item For each vertex $v\in\mc{V}_{h}$, let $e_1,\cdots,e_l$ be the edges in $\mc{E}_h$ having $v$ as an end point, and let $K_1,\cdots,K_l$ be the elements in $\mc{P}_h$ having $v$ as a vertex. For each $i = 1,\cdots,l$, let $\vn_{e_i,v}$ be a unit vector normal to $e_i$ pointing in the counterclockwise direction with respect to the vertex $v$ (see \Cref{fig:ConstructD1}). Define $\vpsi_v\in\vV_h$ by
\begin{dis}
\vpsi_v := h\sum_{i=1}^l\frac{\inn{\vn_{e_i},\vn_{e_i,v}}}{|e_i|}\vvph_{e_i,1}^n + \sum_{j=1}^l\sum_{q\in M_{k-1}(K_j)\setminus\{1\}}c_{j,q}\vvph_{K_j,\nabla q},
\end{dis}
where
\begin{dis}
c_{j,q} = \frac{h}{|K_j|}\sum_{i=1}^l\frac{\inn{\vn_{e_i},\vn_{e_i,v}}}{|e_i|}\int_{e_i\cap \pd K_j}\vvph_{e_i,1}^n\cdot\vn_{K_j}q\diff s
\end{dis}
for $q\in M_{k-1}(K_j)\setminus\{1\}$ and $j = 1,\cdots,l$.
\item For each edge $e\in\mc{E}_{h}$ and each $q\in M_{k-1}(e)$, define $\vpsi_{e,q}^t$ by
\begin{dis}
\vpsi_{e,q}^t = \vvph_{e,q}^t.
\end{dis}
\item Assume $k \geq 2$. For each edge $e\in\mc{E}_{h}$ and each $q\in M_{k-1}(e)\setminus\{1\}$, define $\vpsi_{e,q}^n$ by
\begin{dis}
\vpsi_{e,q}^n = \vvph_{e,q}^n + \sum_{K\in\mc{P}_h}\sum_{r\in M_{k-1}(K)\setminus\{1\}}c_{K,r}\vvph_{K,\nabla r},
\end{dis}
where
\begin{dis}
c_{K,r} = \frac{1}{|K|}\int_{e\cap \pd K}\vvph_{e,q}^n\cdot\vn_Kr\diff s, \quad r\in M_{k-1}(K)\setminus\{1\}, \ K\in\mc{P}_h.
\end{dis}
(See \Cref{fig:ConstructD3}.)
\item Assume $k \geq 3$. For each $K\in\mc{P}_h$ and each $\vq \in \vM_k^\oplus$, define $\vpsi_{K,\vq}$ by
\begin{dis}
\vpsi_{K,\vq} = \vvph_{K,\vq}.
\end{dis}
\end{enumerate}

\begin{remark}
The coefficients $c_{j,q}$ and $c_{K,r}$ defined in (D1) and (D3) are exactly computable using the DOFs \eqref{eqn:EdgeNDOFG}-\eqref{eqn:CellDOFG}. Moreover, since it is computed elementwise, the cost of computing the coefficients is negligible.
\end{remark}

We first show that the functions defined in (D1)-(D4) are indeed contained in $\vZ_h$. 

\begin{lemma}\label{lem:D1D4DivFree}
The functions defined in (D1)-(D4) are contained in $\vZ_h$. 
\end{lemma}

\begin{proof}
Since $\div\vV_h(K)\sus P_{k-1}(K)$ for each $K\in\mc{P}_h$, if $\vv_h\in\vV_h$ then
\begin{eqn}\label{eqn:DivFreeIFF}
\vv_h\in\vZ_h \quad \tr{if and only if} \quad \int_Kq\div\vv_h\diff\vx = 0 \ \forall q\in M_{k-1}(K), \ \forall K\in\mc{P}_h.
\end{eqn}
From \eqref{eqn:DivFreeIFF}, the functions in (D2) and (D4) are obviously contained in $\vZ_h$. We first show that the functions in (D1) belong to $\vZ_h$. Note that
\begin{dis}
\int_Kq\div\vpsi_v\diff\vx = 0 \quad \forall q\in M_{k-1}(K), \ \forall K\in\mc{P}_h \ \tr{with} \ K\neq K_1,\cdots,K_l.
\end{dis}
Let $j = 1,\cdots,l$. Since $K_j$ is a polygon having $v$ as a vertex, there are exactly two edges $e_{i_1,v}$ and $e_{i_2,v}$ with $1\leq i_1,i_2 \leq l$ such that $e_{i_1,v},e_{i_2,v}\sus \pd K_j$. Moreover, one of the normal vectors $\vn_{e_1,v}$ and $\vn_{e_2,v}$ coincides with $\vn_{K_j}$, and the other has the opposite direction of $\vn_{K_j}$. We may assume that $\vn_{e_{i_1,v}} = \vn_{K_j}|_{e_{i_1}}$ and $\vn_{e_{i_2,v}} = - \vn_{K_j}|_{e_{i_2}}$. Then
\begin{eqnarray*}
&& \int_{K_j}\div\vpsi_v\diff\vx = \int_{\pd K_j}\vpsi_v\cdot\vn_{K_j}\diff s \\
& = & h\frac{\inn{\vn_{e_{i_1}},\vn_{e_{i_1,v}}}}{|e_{i_1}|}\int_{e_{i_1}}\vvph_{e_{i_1},1}^n\cdot\vn_{K_j}|_{e_{i_1}}\diff s + h\frac{\inn{\vn_{e_{i_2}},\vn_{e_{i_2,v}}}}{|e_{i_2}|}\int_{e_{i_2}}\vvph_{e_{i_2},1}^n\cdot\vn_{K_j}|_{e_{i_2}}\diff s \\
& = & h\inn{\vn_{e_{i_1}},\vn_{e_{i_1,v}}}\inn{\vn_{K_j}|_{e_{i_1}},\vn_{e_{i_1}}} + h\inn{\vn_{e_{i_2}},\vn_{e_{i_2,v}}}\inn{\vn_{K_j}|_{e_{i_2}},\vn_{e_{i_2}}} \\
& = & h\inn{\vn_{e_{i_1}},\vn_{e_{i_1,v}}}\inn{\vn_{e_{i_1,v}},\vn_{e_{i_1}}} - h\inn{\vn_{e_{i_2}},\vn_{e_{i_2,v}}}\inn{\vn_{e_{i_2,v}},\vn_{e_{i_2}}} \\
& = & 0.
\end{eqnarray*}
Suppose $q\in M_{k-1}(K_j)\setminus\{1\}$. Then
\begin{eqnarray*}
&& \int_{K_j}q\div\vpsi_v\diff\vx \\
& = & \int_{\pd K_j}q\vpsi_v\cdot\vn_K\diff s - \int_{K_j}\vpsi_v\cdot\nabla q\diff\vx \\
& = & h\sum_{i=1}^l\frac{\inn{\vn_{e_i},\vn_{e_i,v}}}{|e_i|}\int_{e_i\cap \pd K_j}q\vvph_{e_i,1}^n\cdot\vn_{K_j}\diff s - \sum_{q'\in M_{k-1}(K_j)\setminus\{1\}}c_{j,q'}\int_{K_j}\vvph_{K_j,\nabla q'}\cdot\nabla q\diff\vx \\
& = & h\sum_{i=1}^l\frac{\inn{\vn_{e_i},\vn_{e_i,v}}}{|e_i|}\int_{e_i\cap \pd K_j}q\vvph_{e_i,1}^n\cdot\vn_{K_j}\diff s - c_{j,q}\int_{K_j}\vvph_{K_j,\nabla q}\cdot\nabla q\diff\vx \\
& = & h\sum_{i=1}^l\frac{\inn{\vn_{e_i},\vn_{e_i,v}}}{|e_i|}\int_{e_i\cap \pd K_j}q\vvph_{e_i,1}^n\cdot\vn_{K_j}\diff s - |K_j|c_{j,q} \\
& = & 0.
\end{eqnarray*}
Here we used the relations
\begin{dis}
\int_{K_j}\vvph_{K_j,\nabla q'}\cdot\nabla q\diff\vx = \left\{\begin{array}{ll}
|K_j| & \tr{if $q = q'$} \\
0 & \tr{if $q \neq q'$}.
\end{array}\right.
\end{dis}
Thus $\vpsi_v\in \vZ_{h}$. We next show that the functions $\vpsi_{e,q}^n$ in (D3) belong to $\vZ_h$. Note that $c_{K,r} = 0$ for any $r\in M_{k-1}(K)\setminus\{1\}$ and any $K\in\mc{P}_h$ with $e\not\sus\pd K$. Then
\begin{eqnarray*}
\int_Kr\div\vpsi_{e,q}^n\diff\vx & = & \int_{\pd K}r\vpsi_{e,q}^n\cdot\vn_K\diff s - \int_K\vpsi_{e,q}^n\cdot\nabla r\diff\vx \\
& = & \int_{\pd K}r\vpsi_{e,q}^n\cdot\vn_K\diff s = 0
\end{eqnarray*}
for any $r\in M_{k-1}(K)\setminus\{1\}$ and $K\in\mc{P}_h$ with $e\not\sus\pd K$. We next suppose that $K\in\mc{P}_h$ satisfies $e\sus\pd K$. Since $q\neq 1$,
\begin{dis}
\int_K\div\vpsi_{e,q}^n\diff\vx = \int_{\pd K}\vpsi_{e,q}^n\cdot\vn_K\diff s = 0,
\end{dis}
and
\begin{eqnarray*}
\int_Kr\div\vpsi_{e,q}^n\diff\vx & = & \int_{\pd K}r\vpsi_{e,q}^n\cdot\vn_K\diff s - \int_K\vpsi_{e,q}^n\cdot\nabla r\diff\vx \\
& = & \int_{e\cap\pd K}r\vvph_{e,q}^n\cdot\vn_K\diff s - \sum_{r'\in M_{k-1}(K)\setminus\{1\}}c_{K,r'}\int_K\vvph_{K,\nabla r'}\cdot\nabla r\diff\vx \\
& = & \int_{e\cap\pd K}r\vvph_{e,q}^n\cdot\vn_K\diff s - c_{K,r}\int_K\vvph_{K,\nabla r}\cdot\nabla r\diff\vx \\
& = & \int_{e\cap\pd K}r\vvph_{e,q}^n\cdot\vn_K\diff s - |K|c_{K,r} \\
& = & 0
\end{eqnarray*}
for any $r\in M_{k-1}(K)\setminus \{1\}$. Here, as before, we used the relations
\begin{dis}
\int_K\vvph_{K,\nabla r'}\cdot\nabla r\diff\vx = \left\{\begin{array}{ll}
|K| & \tr{if $r = r'$} \\
0 & \tr{if $r \neq r'$}.
\end{array}\right.
\end{dis}
Thus $\vpsi_{e,q}^n\in\vZ_{h}$. This concludes the proof of the lemma.
\end{proof}

The next theorem shows that some of these functions generate a basis for $\vZ_{h,0}$. 

\begin{theorem}\label{thm:ConstructDivFree}
Let $\vZ_1,\vZ_2,\vZ_3,\vZ_4$ be the subspaces of $\vV_{h}$ defined by
\begin{eqnarray*}
\vZ_1 & = & \spn\left(\left\{\vpsi_v : v\in\mc{V}_h^i\right\}\right), \\
\vZ_2 & = & \spn\left(\left\{\vpsi_{e,q}^t : q\in M_{k-1}(e), e\in\mc{E}_h^i\right\}\right), \\
\vZ_3 & = & \left\{\begin{array}{ll}
\spn\left(\left\{\vpsi_{e,q}^n : q\in M_{k-1}(e)\setminus\{1\}, e\in\mc{E}_h^i\right\}\right) & \tr{if $k \geq 2$} \\
\{0\} & \tr{otherwise}
\end{array}\right.,\\
\vZ_4 & = & \left\{\begin{array}{ll}
\spn\left(\left\{\vpsi_{K,\vq} : \vq\in \vM_k^\oplus, K\in\mc{P}_h\right\}\right) & \tr{if $k \geq 3$} \\
\{0\} & \tr{otherwise}
\end{array}\right.,
\end{eqnarray*}
where $\vpsi_v$, $\vpsi_{e,q}^t$, $\vpsi_{e,q}^n$, and $\vpsi_{K,\vq}$ are the functions given in (D1)-(D4), respevtively. Then the following hold.
\begin{enumerate}[(i)]
\item $\vZ_1,\vZ_2,\vZ_3,\vZ_4\sus\vZ_{h,0}$.
\item $\vZ_i\cap \vZ_j = \{0\}$ for any pair $(i,j)$ with $i\neq j$.
\item The dimensions of the subspaces $\vZ_1$, $\vZ_2$, $\vZ_3$, and $\vZ_4$ satisfy
\begin{dis}
\dim \vZ_1 = N_{V,i}, \ \dim \vZ_2 = kN_{E,i}, \ \dim \vZ_3 = (k-1)N_{E,i}, \ \dim \vZ_4 = \frac{(k-1)(k-2)}{2}N_P.
\end{dis}
\end{enumerate}
Consequently, $\vZ_{h,0} = \vZ_1 \oplus \vZ_2 \oplus \vZ_3 \oplus \vZ_4$. 
\end{theorem}

\begin{proof}
Since $\vZ_{h,0} = \vZ_h\cap\vV_{h,0}$, and from \Cref{lem:D1D4DivFree}, it suffices to show that the functions in (D1)-(D4) are contained in $\vV_{h,0}$. Clearly $\vpsi_{K,\vq}\in\vV_{h,0}$ for any $K\in\mc{P}_h$ and any $\vq\in \vM_k^\oplus$. If $e\in\mc{E}_h^i$, then $\vpsi_{e,q}^t\in\vV_{h,0}$ for any $q\in M_{k-1}(e)$ and $\vpsi_{e,q}^n\in\vV_{h,0}$ for any $q\in M_{k-1}(e)\setminus\{1\}$. If $v\in\mc{V}_h^i$, then the edges in $\mc{E}_h$ that have $v$ as an end point are contained in $\mc{E}_h^i$. Thus $\vpsi_v\in\vV_{h,0}$ for any $v\in\mc{V}_h^i$. Hence $\vZ_1$, $\vZ_2$, $\vZ_3$, and $\vZ_4$ are subspaces of $\vZ_{h,0}$.

On the other hand, it is easy to show that $\vZ_i\cap \vZ_j = \{0\}$ for any pair $(i,j)$ with $i\neq j$ and
\begin{dis}
\dim \vZ_1 = N_{V,i}, \ \dim \vZ_2 = kN_{E,i}, \ \dim \vZ_3 = (k-1)N_{E,i}, \ \dim \vZ_4 = \frac{(k-1)(k-2)}{2}N_P.
\end{dis}
Then, since $\vZ_1 \oplus \vZ_2 \oplus \vZ_3 \oplus \vZ_4 \sus \vZ_{h,0}$ and since $\dim\vZ_{h,0} = \dim\vZ_1 + \dim\vZ_2 + \dim\vZ_3 + \dim\vZ_4$ by \Cref{prop:DimZh}, we obtain
\begin{dis}
\vZ_{h,0} = \vZ_1 \oplus \vZ_2 \oplus \vZ_3 \oplus \vZ_4.
\end{dis}
This concludes the proof of the theorem.
\end{proof}

\begin{remark}
If $k = 1$ and the mesh $\mc{P}_h$ is a triangular mesh, then the construction of the basis of $\vZ_{h,0}$ described in \Cref{thm:ConstructDivFree} is exactly the same with the divergence-free basis in the Crouzeix-Raviart finite element space \cite{thomassetimplementation,brenner1990nonconforming}.
\end{remark}

\section{Implementation Details}

In this section, we present how to compute the solution $(\vu_h,p_h)$ of the discrete problem \eqref{eqn:StokesDis} by using the construction of $\vZ_{h,0}$ presented in \Cref{sec:NCVEM1}. 

\subsection{Computing the discrete velocity $\vu_h$}

We first consider the case $\vg = \zz$. Note that the discrete velocity $\vu_h$ is the solution of the following discrete problem \cite{crouzeix1973conforming}: Find $\vu_h\in\vZ_{h,0}$ such that
\begin{eqn}\label{eqn:StokesDisVel}
a_h(\vu_h,\vv_h) = \inn{\vf_h,\vv_h} \quad \forall \vv_h\in \vZ_{h,0}.
\end{eqn}
Since $\dim\vZ_{h,0} = \dim\vV_{h,0} - \dim Q_h$, the system \eqref{eqn:StokesDisVel} has a smaller number of unknowns than system \eqref{eqn:StokesDis}. Moreover, system \eqref{eqn:StokesDisVel} is symmetric positive definite, while problem \eqref{eqn:StokesDis} is a saddle point problem. Thus, it is more efficient to compute $\vu_h$ from \eqref{eqn:StokesDisVel} than from the problem \eqref{eqn:StokesDis}.

We next consider the case $\vg \neq \zz$. Let us decompose $\vu_h\in\vV_{h,\vg}$ into
\begin{dis}
\vu_h = \vu_{h,0} + \wt{\vu}_{h},
\end{dis}
where $\wt{\vu}_h\in \vV_{h,\vg} \cap \vZ_h$ and $\vu_{h,0}\in \vZ_{h,0}$ is the solution of the problem
\begin{dis}
a_h(\vu_{h,0},\vv_h) = \inn{\vf_h,\vv_h} - a_h(\wt{\vu}_h,\vv_h) \quad \forall \vv_h\in \vZ_{h,0}.
\end{dis}
Using the construction of $\vZ_{h,0}$ presented in \Cref{thm:ConstructDivFree}, we can compute $\vu_{h,0}$ by solving a symmetric positive definite system of linear equations, as explained in the case $\vg = \zz$. It remains to find a function $\wt{\vu}_h\in \vV_{h,\vg} \cap \vZ_h$. The following theorem shows that we can easily find such a function.

\begin{theorem}\label{thm:DivFreeNonzeroBdry}
Let $N = N_V^\pd$ and label the vertices in $\mc{V}_h^{\pd}$ by $1,2,\cdots,N$ such that $v_1,\cdots,v_N$ are in counterclockwise order with respect to $\Omega$. We also label the edges in $\mc{E}_h^{\pd}$ by $1,2,\cdots,N$, such that the endpoints of the edge $e_i$ are $v_i$ and $v_{i+1}$ for $i = 1,2,\cdots,N-1$, and the endpoints of the edge $e_{N}$ are $v_{N}$ and $v_1$ (since $\Omega$ is a simply connected polygon, $N_V^{\pd} = N_E^{\pd}$). Let $\wt{\vu}_h$ be the function in $\vV_h$ defined by
\begin{dis}
\wt{\vu}_h = \sum_{v\in\mc{V}_h^{\pd}}C_{1,v}\vpsi_{v} + \sum_{e\in\mc{E}_h^{\pd}}\sum_{q\in M_{k-1}(e)}C_{2,e,q}\vpsi_{e,q}^t + \sum_{e\in\mc{E}_h^{\pd}}\sum_{q\in M_{k-1}(e)\setminus\{1\}}C_{3,e,q}\vpsi_{e,q}^n,
\end{dis}
where the coefficients $(C_{1,v})_v = (C_{1,v_1},\cdots,C_{1,v_N})$ are given by
\begin{dis}
C_{1,v_k} = -\sum_{i=k}^{N}\int_{e_i}\vg\cdot\vn_{e_i}\diff s, \quad k = 1,2,\cdots,N,
\end{dis}
and the coefficients $(C_{2,e,q})_{e,q}$, and $(C_{3,e,q})_{e,q}$ are given by
\begin{eqnarray*}
C_{2,e,q} & = & \frac{1}{|e|}\int_e\vg\cdot\vt_eq\diff s, \quad q\in M_{k-1}(e), \ e\in\mc{E}_h^{\pd}, \\
C_{3,e,q} & = & \frac{1}{|e|}\int_e\vg\cdot\vn_eq\diff s, \quad q\in M_{k-1}(e)\setminus\{1\}, \ e\in\mc{E}_h^{\pd}.
\end{eqnarray*}
Then $\wt{\vu}_h\in \vV_{h,\vg}\cap \vZ_h$. 
\end{theorem}

\begin{proof}
From the construction of $\wt{\vu}_h$, it is obvious that $\div_h\wt{\vu}_h = 0$. Thus it remains to show that $\wt{\vu}_h\in \vV_{h,\vg}$. From the definition of the coefficients $(C_{2,e,q})_{e,q}$, and $(C_{3,e,q})_{e,q}$, we obtain
\begin{eqnarray*}
\int_e\wt{\vu}_h\cdot\vt_eq\diff s & = & \int_e\vg\cdot\vt_eq\diff s, \quad q\in M_{k-1}(e), \ e\in\mc{E}_h^{\pd}, \\
\int_e\wt{\vu}_h\cdot\vn_eq\diff s & = & \int_e\vg\cdot\vn_eq\diff s, \quad q\in M_{k-1}(e)\setminus\{1\}, \ e\in\mc{E}_h^{\pd}.
\end{eqnarray*}
Since the boundary edge $e_i$ with $1\leq i \leq N-1$ has endpoints $v_i$ and $v_{i+1}$, and since the vertices $v_1,\cdots,v_N$ are labeled in counterclockwise order with respect to $\Omega$, we obtain
\begin{dis}
\vn_{e_i,v_{i+1}} = \vn_{e_i} = -\vn_{e_i,v_{i}},
\end{dis}
where $\vn_{e_i}$ is a unit normal vector in the outward direction with respect to $\Omega$, and $\vn_{e_i,v_i}$ and $\vn_{e_i,v_{i+1}}$ are unit vectors normal to $e_i$ pointing in the counterclockwise direction with respect to $v_i$ and $v_{i+1}$, respectively (see \Cref{fig:NormalBdry}). Similarly, we obtain
\begin{dis}
\vn_{e_N,v_1} = \vn_{e_N} = -\vn_{e_N,v_N}.
\end{dis}
Thus we have
\begin{eqnarray*}
\int_{e_i}\wt{\vu}_h\cdot\vn_{e_i}\diff s & = & -C_{1,v_i} + C_{1,v_{i+1}} \quad \forall i = 1,2,\cdots,N-1, \\
\int_{e_N}\wt{\vu}_h\cdot\vn_{e_N}\diff s & = & -C_{1,v_{N}} + C_{1,v_1}.
\end{eqnarray*}
Using the definition of the coefficients $(C_{1,v})$, 
\begin{dis}
-C_{1,v_i} + C_{1,v_{i+1}} = \int_{e_i}\vg\cdot\vn_{e_i}\diff s \quad \forall i = 1,2,\cdots,N-1. 
\end{dis}
From \eqref{eqn:StokesBdry} we obtain $C_{1,v_1} = -\int_{\pd\Omega}\vg\cdot\vn_{\Omega}\diff s = 0$ and thus
\begin{eqnarray*}
-C_{1,v_{N}} + C_{1,v_{1}} = \int_{e_{N}}\vg\cdot\vn_{e_{N}}\diff s.
\end{eqnarray*}
Therefore $\wt{\vu}_h\in \vV_{h,\vg}$. 
\end{proof}

\begin{figure}
\begin{center}
\begin{tikzpicture}
\fill (-1,-1.5) circle [radius = 0.075];
\fill (1.5,0) circle [radius = 0.075];
\draw [thick] (-3,-1.5) -- (-1,-1.5);
\draw [thick] (-1,-1.5) -- (1.5,0);
\draw [thick, dashed] (1.5,0) -- (0.5,1.5);
\draw [thick] (1.5,0) -- (3,1.75);
\draw [thick, dashed] (-1,-1.5) -- (-1.5,0.5);
\draw [thick, dashed] (0.5,1.5) -- (0.5,2.5);
\draw [thick, dashed] (-1.5,0.5) -- (-3,1);
\draw [thick, dashed] (-1.5,0.5) -- (-1,1.5);
\draw [thick, dashed] (-1,1.5) -- (0.5,1.5);
\draw [thick, ->] (0.25,-0.75) -- (1,-2);
\draw [thick, ->] (0.25,-0.75) -- (-0.5,0.5);
\node at (-2.75,2.25) {$\Omega$};
\node at (2.75,-2.25) {$\Omega^c$};
\node at (-1,-1.75){$v_{i}$};
\node[right] at (1.5,-0.25){$v_{i+1}$};
\node[right] at (0.8,-1.25){$\vn_{e_i} = \vn_{e_i,v_{i+1}}$};
\node[right] at (-0.25,0.25){$\vn_{e_i,v_{i}}$};
\node at (0.375,-0.5){$e_i$};
\end{tikzpicture}
\caption{}
\label{fig:NormalBdry}
\end{center}
\end{figure}

\subsection{Recovery of the discrete pressure $p_h$}

Once we have the discrete velocity $\vu_h$, the discrete pressure $p_h$ can be obtained by solving the overdetermined system
\begin{eqn}\label{eqn:StokesDisPre}
b_h(\vv_h,p_h) = \inn{\vf_h,\vv_h} - a_h(\vu_h,\vv_h) \quad \forall \vv_h\in\vV_{h,0}.
\end{eqn}

\section{Numerical Experiments}

In this section, we present several numerical experiments for the symmetric positive definite linear system \eqref{eqn:StokesDisVel} and the overdetermined linear system \eqref{eqn:StokesDisPre}. Consider the Stokes problem \eqref{eqn:Stokes123} on the unit square domain $\Omega = [0,1]^2$, where the exact solution is given by
\begin{eqnarray*}
\vu(x,y) & = & ((1 - \cos(2\pi x))\sin(2\pi y), -(1-\cos(2\pi y))\sin(2\pi x)), \\
p(x,y) & = & e^x - e^y.
\end{eqnarray*}
We solve both \eqref{eqn:StokesDisVel} and \eqref{eqn:StokesDisPre} for $k = 1,2,3$, and we compute the velocity error in the discrete energy norm
\begin{dis}
E_v := a_h(\vu_h - I_h\vu,\vu_h - I_h\vu)^{1/2}
\end{dis}
and the pressure error in the $L^2$-norm
\begin{dis}
E_p := \|p_h - \Pi_hp\|_{0,\Omega},
\end{dis}
where $\Pi_hp$ is the piecewise polynomial function such that for each $K\in\mc{P}_h$ the restriction $\Pi_hp|_K$ is the $L^2$-projection of $p$ onto $P_{k-1}(K)$. 

\begin{figure}[!ht]
\centering
\includegraphics[width = 0.33\textwidth]{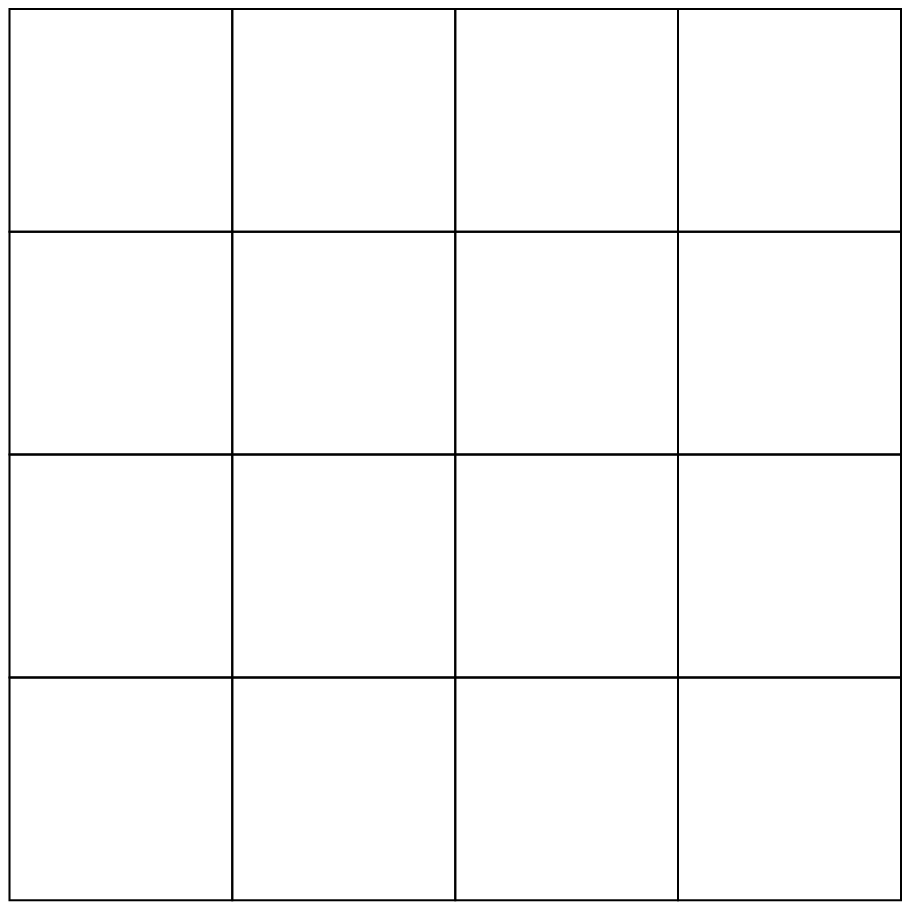}
\includegraphics[width = 0.33\textwidth]{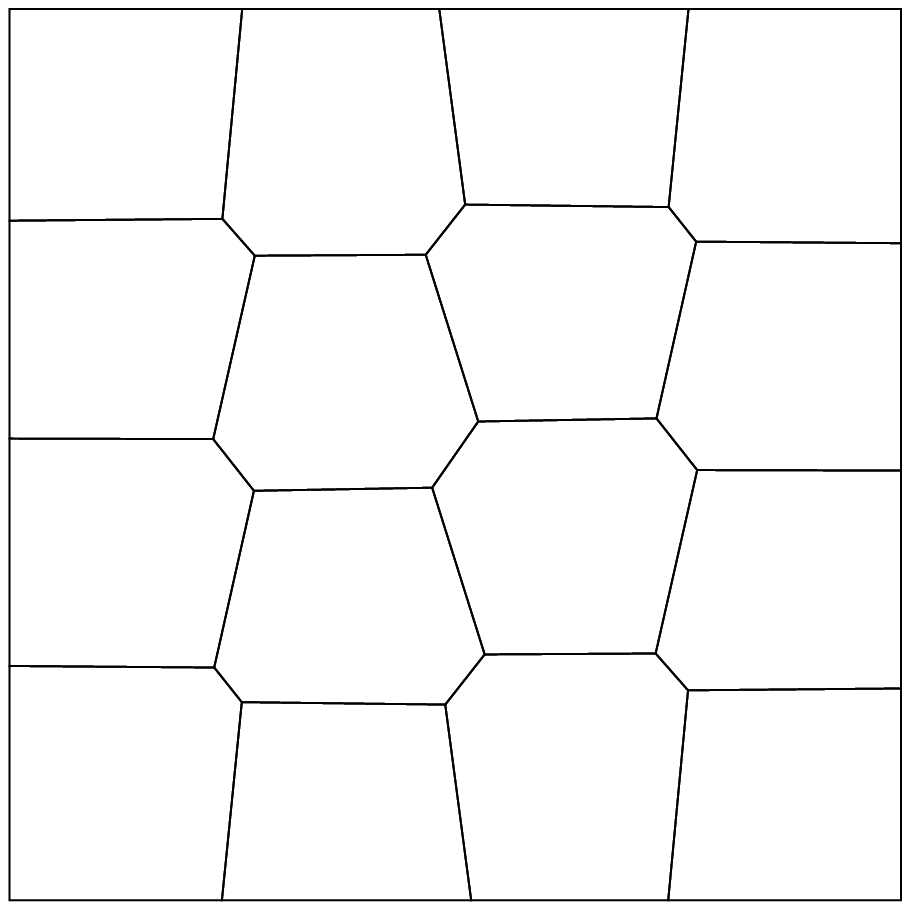}
\includegraphics[width = 0.33\textwidth]{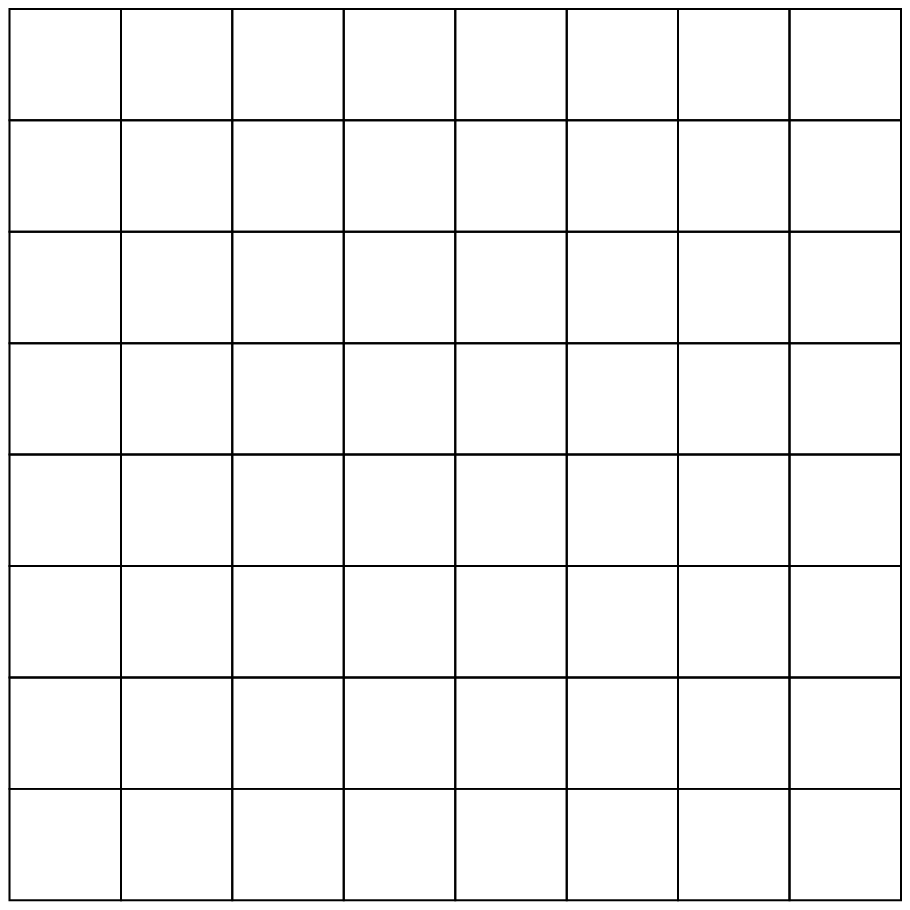}
\includegraphics[width = 0.33\textwidth]{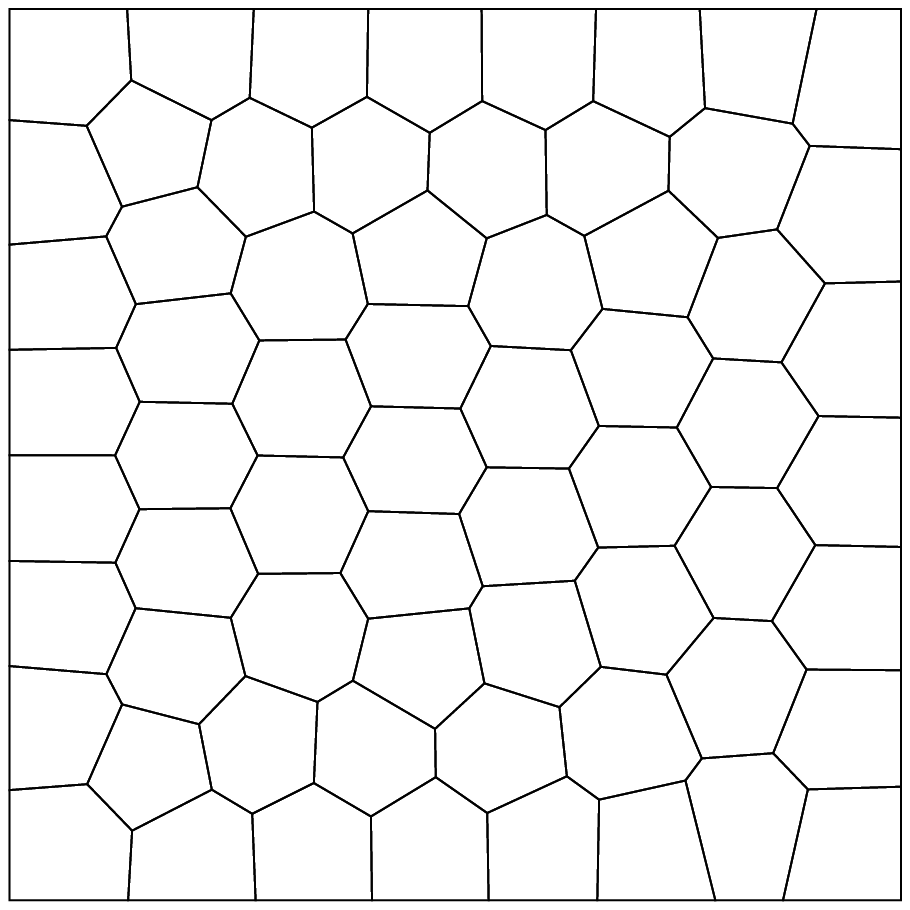}
\caption{The meshes $\mc{P}_h^1$ (left), and $\mc{P}_h^2$ (right).}
\label{fig:mesh}
\end{figure}

We decompose $\Omega$ into the following sequences of convex polygonal meshes:
\begin{enumerate}[(i)]
\item uniform square meshes $\mc{P}_h^1$ with $h = 1/4$, $1/8$, $1/16$, $1/32$, $1/64$, $1/128$,
\item unstructured polygonal meshes $\mc{P}_h^2$ with $h = 1/4$, $1/8$, $1/16$, $1/32$, $1/64$, $1/128$.
\end{enumerate}
Some examples of the meshes are shown in \Cref{fig:mesh}. The unstructured polygonal meshes $\{\mc{P}_h^2\}_h$ are generated from PolyMesher \cite{talischi2012polymesher}. Mesh data (the number of polygons, interior edges, and interior vertices) for each $h$ are given in \Cref{tab:MeshInfo}. 

\begin{table}[!ht]
\footnotesize
\centering
\caption{Mesh information}
\label{tab:MeshInfo}
\begin{tabular}{|c||c|c|c||c|c|c|}
\hline
\multirow{2}{*}{$h$} & \multicolumn{3}{c||}{$\mc{P}_h^1$}& \multicolumn{3}{c|}{$\mc{P}_h^2$} \\ \cline{2-7} 
                     & $N_P$   & $N_{E,i}$  & $N_{V,i}$  & $N_P$   & $N_{E,i}$  & $N_{V,i}$  \\ \hline\hline
1/4                  & 16      & 24         & 9          & 16      & 33         & 18         \\ \hline
1/8                  & 64      & 112        & 49         & 64      & 162        & 99         \\ \hline
1/16                 & 256     & 480        & 225        & 256     & 707        & 452        \\ \hline
1/32                 & 1024    & 1984       & 961        & 1024    & 2953       & 1930       \\ \hline
1/64                 & 4096    & 8064       & 3969       & 4096    & 12043      & 7948       \\ \hline
1/128                & 16384   & 32512      & 16129      & 16384   & 48655      & 32272      \\ \hline
\end{tabular}
\end{table}

In \Crefrange{tab:DimDiscreteSpace1}{tab:DimDiscreteSpace3}, we present the dimensions of the spaces $\vV_{h,0}$, $Q_h$, and $\vZ_{h,0}$, for each mesh $\mc{P}_h^1$, $\mc{P}_h^2$ and each $k = 1,2,3$. Since the number of unknowns of the system \eqref{eqn:StokesDis} is $\dim\vV_{h,0} + \dim Q_h$ and the number of unknowns of the system \eqref{eqn:StokesDisVel} is $\dim\vZ_{h,0}$, we can see that the system \eqref{eqn:StokesDisVel} has fewer unknowns than the system \eqref{eqn:StokesDis}.

\begin{table}[!ht]
\footnotesize
\centering
\caption{Dimensions of the discrete spaces ($k = 1$)}
\label{tab:DimDiscreteSpace1}
\begin{tabular}{|c||c|c|c||c|c|c|}
\hline
\multirow{2}{*}{$h$} & \multicolumn{3}{c||}{$\mc{P}_h^1$}              & \multicolumn{3}{c|}{$\mc{P}_h^2$}              \\ \cline{2-7} 
                     & $\dim\vV_{h,0}$ & $\dim Q_h$ & $\dim\vZ_{h,0}$ & $\dim\vV_{h,0}$ & $\dim Q_h$ & $\dim\vZ_{h,0}$ \\ \hline\hline
1/4                  & 48              & 15         & 33              & 66              & 15         & 51              \\ \hline
1/8                  & 224             & 63         & 161             & 324             & 63         & 261             \\ \hline
1/16                 & 960             & 255        & 705             & 1414            & 255        & 1159            \\ \hline
1/32                 & 3968            & 1023       & 2945            & 5906            & 1023       & 4883            \\ \hline
1/64                 & 16128           & 4095       & 12033           & 24086           & 4095       & 19991           \\ \hline
1/128                & 65024           & 16383      & 48641           & 97310           & 16383    & 80927           \\ \hline
\end{tabular}
\end{table}

\begin{table}[!ht]
\footnotesize
\centering
\caption{Dimensions of the discrete spaces ($k = 2$)}
\label{tab:DimDiscreteSpace2}
\begin{tabular}{|c||c|c|c||c|c|c|}
\hline
\multirow{2}{*}{$h$} & \multicolumn{3}{c||}{$\mc{P}_h^1$}              & \multicolumn{3}{c|}{$\mc{P}_h^2$}              \\ \cline{2-7} 
                     & $\dim\vV_{h,0}$ & $\dim Q_h$ & $\dim\vZ_{h,0}$ & $\dim\vV_{h,0}$ & $\dim Q_h$ & $\dim\vZ_{h,0}$ \\ \hline\hline
1/4                  & 128             & 47         & 81              & 164             & 47         & 117             \\ \hline
1/8                  & 576             & 191        & 385             & 776             & 191        & 585             \\ \hline
1/16                 & 2432            & 767        & 1665            & 3340            & 767        & 2573            \\ \hline
1/32                 & 9984            & 3071       & 6913            & 13860           & 3071       & 10789           \\ \hline
1/64                 & 40448           & 12287      & 28161           & 56364           & 12287      & 44077           \\ \hline
1/128                & 162816          & 49151      & 113665          & 227388          & 49151		& 178237          \\ \hline
\end{tabular}
\end{table}

\begin{table}[!ht]
\footnotesize
\centering
\caption{Dimensions of the discrete spaces ($k = 3$)}
\label{tab:DimDiscreteSpace3}
\begin{tabular}{|c||c|c|c||c|c|c|}
\hline
\multirow{2}{*}{$h$} & \multicolumn{3}{c||}{$\mc{P}_h^1$}              & \multicolumn{3}{c|}{$\mc{P}_h^2$}              \\ \cline{2-7} 
                     & $\dim\vV_{h,0}$ & $\dim Q_h$ & $\dim\vZ_{h,0}$ & $\dim\vV_{h,0}$ & $\dim Q_h$ & $\dim\vZ_{h,0}$ \\ \hline\hline
1/4                  & 240             & 95         & 145             & 294             & 95         & 199             \\ \hline
1/8                  & 1056            & 383        & 673             & 1356            & 383        & 973             \\ \hline
1/16                 & 4416            & 1535       & 2881            & 5778            & 1535       & 4243            \\ \hline
1/32                 & 18048           & 6143       & 11905           & 23862           & 6143       & 17719           \\ \hline
1/64                 & 72960           & 24575      & 48385           & 96834           & 24575      & 72259           \\ \hline
1/128                & 293376          & 98303      & 195073          & 390234          & 98303		& 291931          \\ \hline
\end{tabular}
\end{table}

\begin{figure}[!ht]
\centering
\includegraphics[width = 0.4\textwidth]{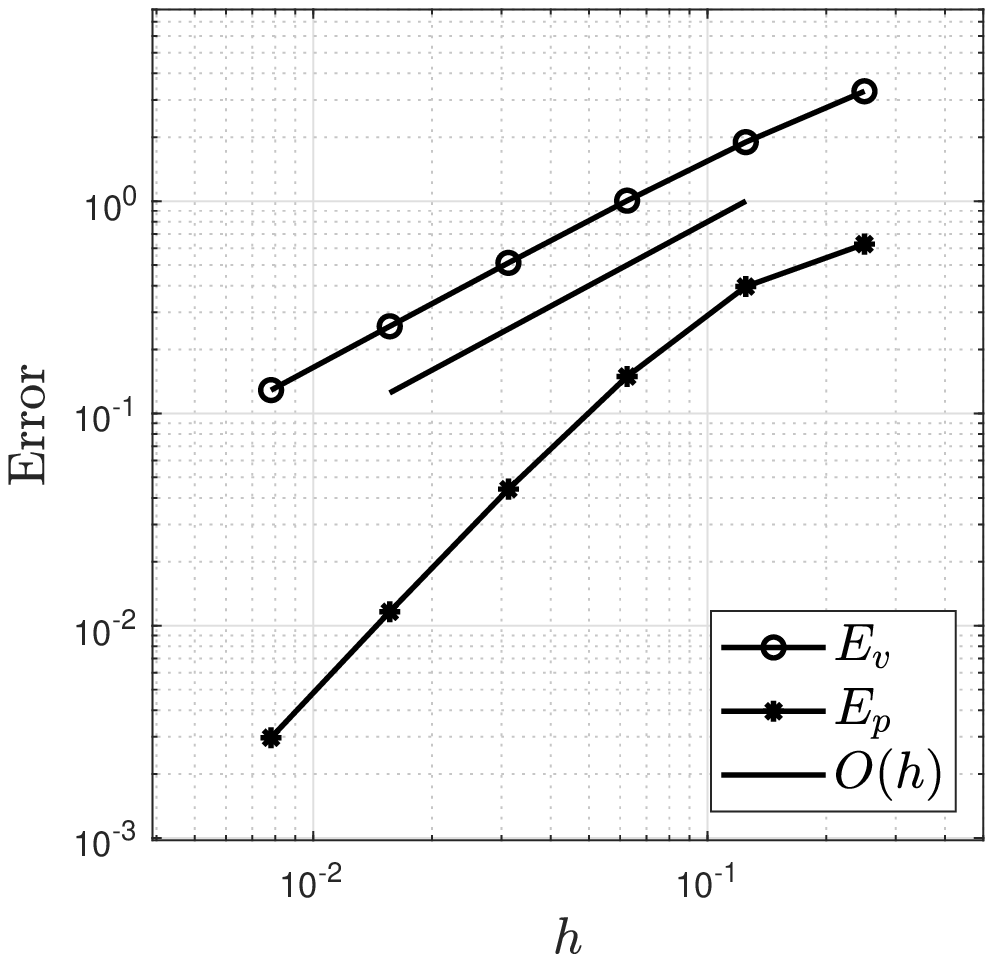}
\includegraphics[width = 0.4\textwidth]{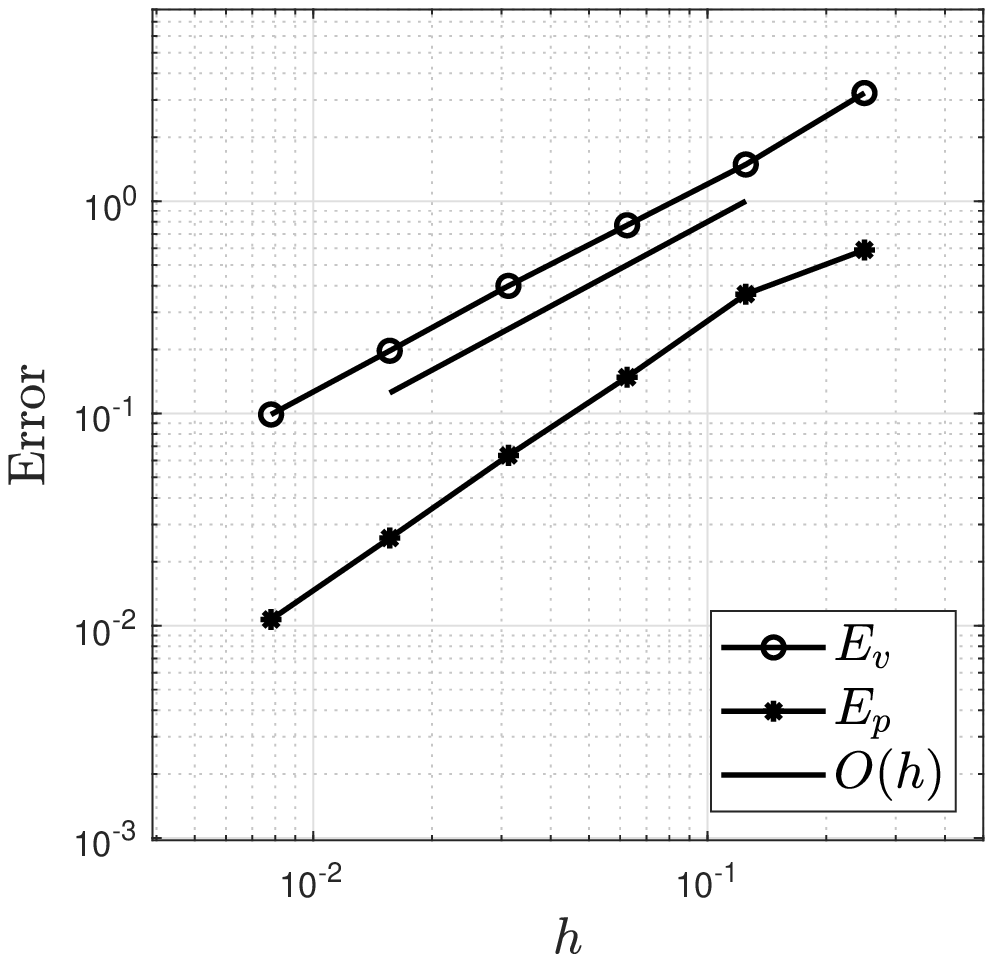}
\caption{Error curves with respect to $h$ for the velocity and pressure on the sequences of meshes $\mc{P}_h^1$ (left) and $\mc{P}_h^2$ (right) with $k = 1$.}
\label{fig:error1}
\end{figure}

\begin{figure}[!ht]
\centering
\includegraphics[width = 0.4\textwidth]{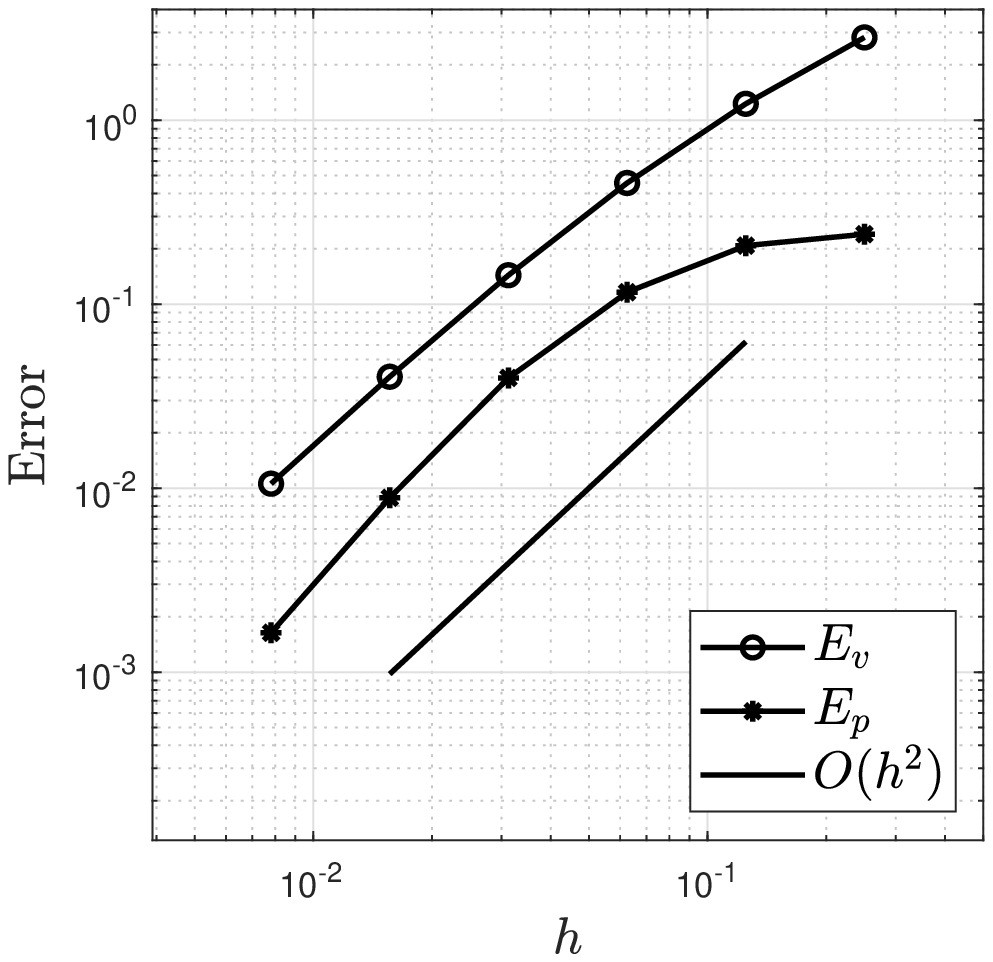}
\includegraphics[width = 0.4\textwidth]{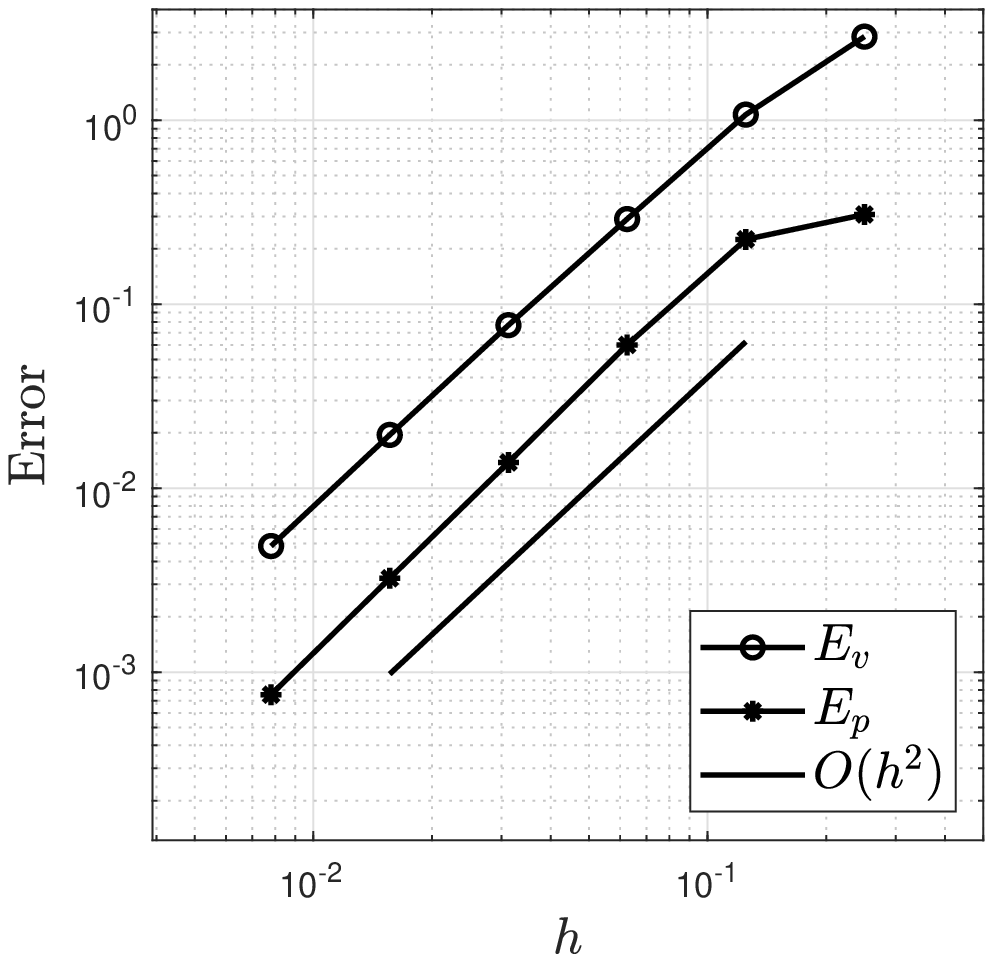}
\caption{Error curves with respect to $h$ for the velocity and pressure on the sequences of meshes $\mc{P}_h^1$ (left) and $\mc{P}_h^2$ (right) with $k = 2$.}
\label{fig:error2}
\end{figure}

\begin{figure}[!ht]
\centering
\includegraphics[width = 0.4\textwidth]{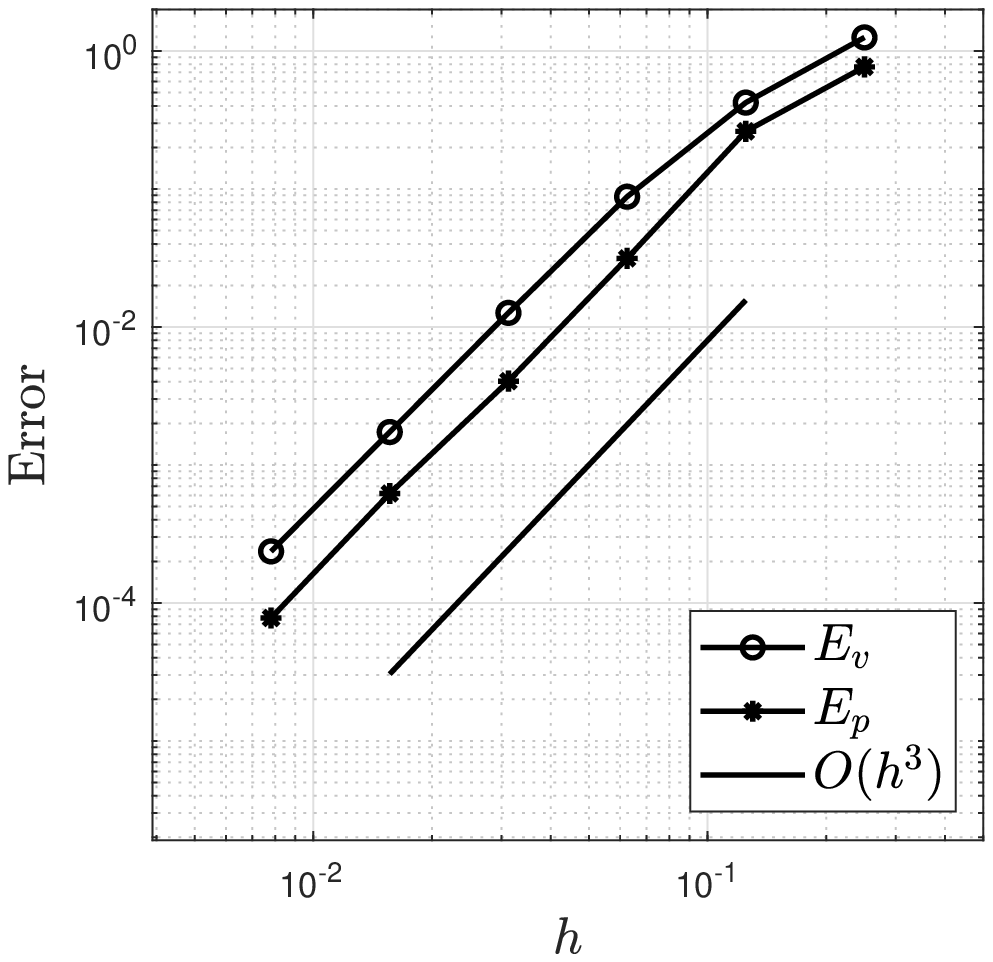}
\includegraphics[width = 0.4\textwidth]{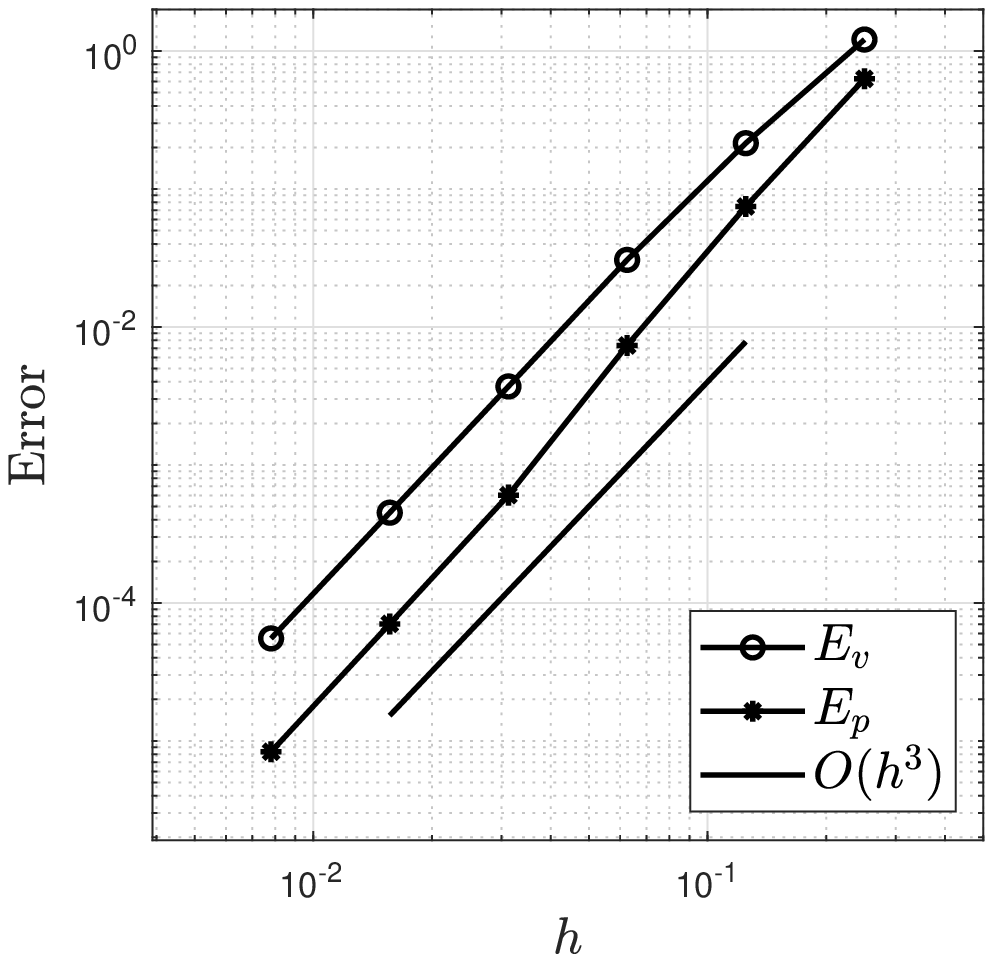}
\caption{Error curves with respect to $h$ for the velocity and pressure on the sequences of meshes $\mc{P}_h^1$ (left) and $\mc{P}_h^2$ (right) with $k = 3$.}
\label{fig:error3}
\end{figure}

The errors $E_v$ and $E_p$ and their orders on the sequences of the meshes for $k = 1,2,3$ are given in \Crefrange{fig:error1}{fig:error3}. In these figures, we see that the convergence order of the errors $E_v$ and $E_p$ are
$O(h^k)$ for $k = 1,2,3$. Thus the numerical results confirm the theoretical analysis in \Cref{thm:Error}. 

In \Cref{tab:CPUTimeTotal}, we compare the CPU running times (on a PC with an Intel Core i5 processor and 8GB RAM) required to solve the reduced system \eqref{eqn:StokesDisVel} and the original saddle-point system \eqref{eqn:StokesDis}, for the uniform square meshes $\{\mc{P}_h^1\}_h$ and $k = 1,2,3$. For a fair comparison, we use unpreconditioned conjugate gradient method (CG) to solve \eqref{eqn:StokesDisVel} and the standard Uzawa method to solve \eqref{eqn:StokesDis}. The cost for computing the discrete pressure (by solving \eqref{eqn:StokesDisPre}) is a fraction of CG, hence it is not included. For each experiment, we write ``$*$'' if the CPU time is more than 518,400 seconds (6 days). For all cases, the CPU time of solving the reduced system is much smaller than that of solving the saddle-point system. 

\begin{table}[!ht]
\footnotesize
\caption{CPU running times}
\label{tab:CPUTimeTotal}
\begin{tabular}{|c||c|c||c|c||c|c|}
\hline
\multirow{3}{*}{$h$} & \multicolumn{6}{c|}{CPU time (secs)}                                                      \\ \cline{2-7}
                     & \multicolumn{2}{c||}{$k=1$}  & \multicolumn{2}{c||}{$k=2$}   & \multicolumn{2}{c|}{$k=3$}   \\ \cline{2-7} 
                     & CG     & Uzawa     & CG      & Uzawa     & CG      & Uzawa     \\ \hline\hline
$1/4$                & 0.001  & 0.095     & 0.001   & 3.274     & 0.004   & 77.824    \\ \hline
$1/8$                & 0.005  & 2.202     & 0.011   & 78.528    & 0.034   & 1713.996  \\ \hline
$1/16$               & 0.034  & 61.202    & 0.066   & 1756.959  & 0.389   & 29664.469 \\ \hline
$1/32$               & 0.275  & 1822.032  & 0.776   & 36869.009 & 4.255   & $*$       \\ \hline
$1/64$               & 4.700  & 53004.713 & 9.581   & $*$       & 53.712  & $*$       \\ \hline
$1/128$              & 70.250 & $*$       & 128.755 & $*$       & 721.477 & $*$       \\ \hline
\end{tabular}
\end{table}

\begin{figure}[!ht]
\centering
\includegraphics[width = 0.4\textwidth]{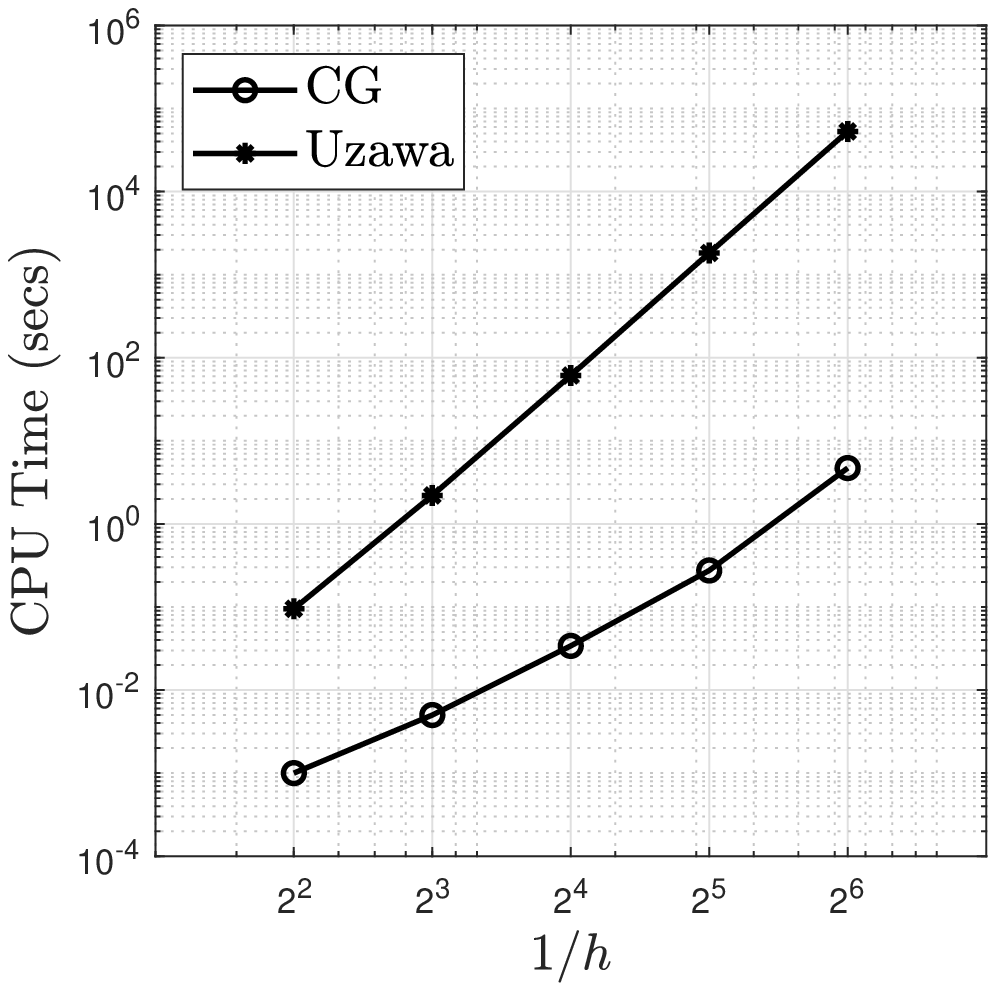}
\includegraphics[width = 0.4\textwidth]{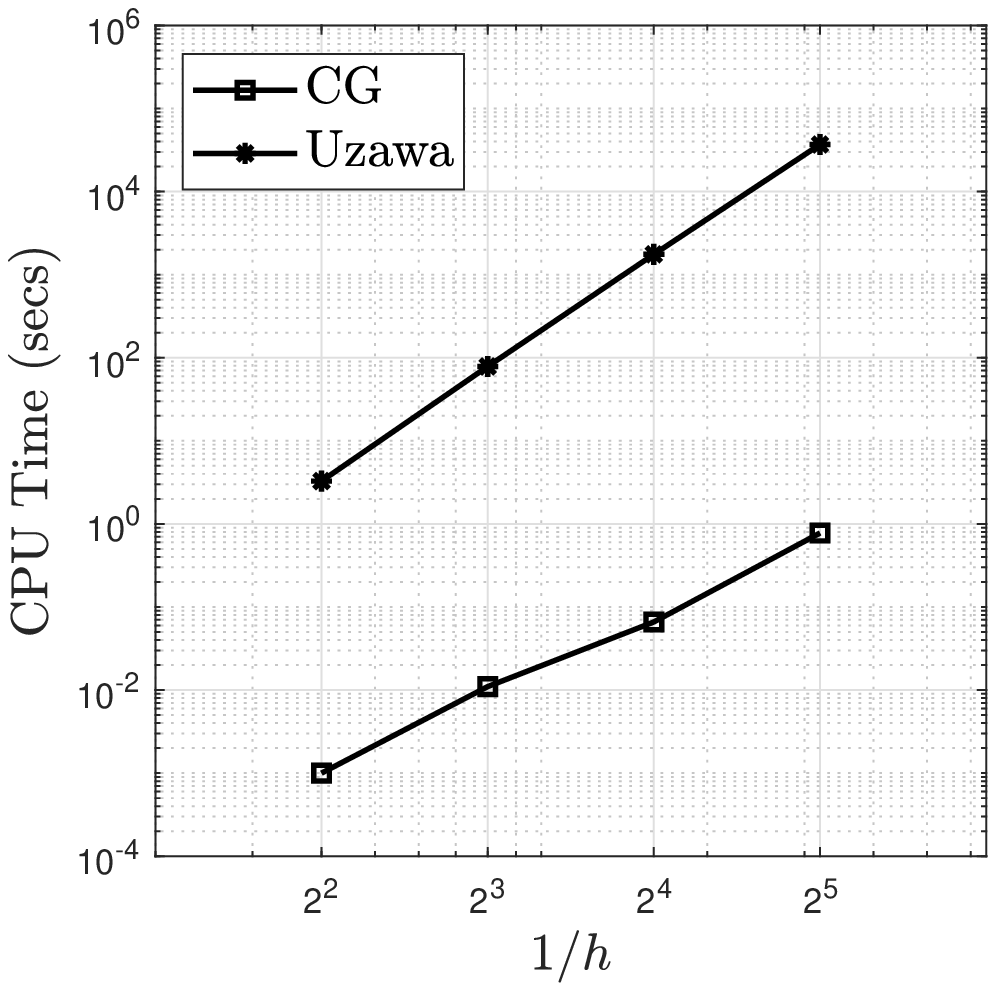}
\caption{CPU time curves of CG and Uzawa with respect to $1/h$ with $k = 1$ (left) and $k = 2$ (right).}
\label{fig:CPUtime}
\end{figure}

\section{Conclusions}

We presented a formal construction of divergence-free bases in the nonconforming VEM for solving the stationary Stokes problem on arbitrary polygonal meshes introduced in \cite{zhao2019divergence}. If $k = 1$ and the mesh is triangular, then the proposed construction of the basis is exactly the same as the divergence-free basis in the Crouzeix-Raviart finite element space \cite{thomassetimplementation,brenner1990nonconforming}. Using our construction, we are able to eliminate the pressure variable from the discrete saddle point formulation, and reduce it to a symmetric positive definite linear system in the velocity variable only. Thus, we can apply many efficient solvers available for symmetric positive definite systems. Finally, we provided some numerical experiments confirming the theoretical results and the efficiency of our construction of divergence-free bases in the nonconforming VEM for the Stokes problem. 

\bibliographystyle{siam}
\bibliography{Refs}

\begin{thebibliography}{10}

\bibitem{ahmad2013equivalent}
{\sc B.~Ahmad, A.~Alsaedi, F.~Brezzi, L.~D. Marini, and A.~Russo}, {\em
  Equivalent projectors for virtual element methods}, Comput. Math. Appl., 66
  (2013), pp.~376--391.

\bibitem{antonietti2014stream}
{\sc P.~F. Antonietti, L.~Beir\~{a}o~da Veiga, D.~Mora, and M.~Verani}, {\em A
  stream virtual element formulation of the {S}tokes problem on polygonal
  meshes}, SIAM J. Numer. Anal., 52 (2014), pp.~386--404.

\bibitem{antonietti2018fully}
{\sc P.~F. Antonietti, G.~Manzini, and M.~Verani}, {\em The fully nonconforming
  virtual element method for biharmonic problems}, Math. Models Methods Appl.
  Sci., 28 (2018), pp.~387--407.

\bibitem{de2016nonconforming}
{\sc B.~Ayuso~de Dios, K.~Lipnikov, and G.~Manzini}, {\em The nonconforming
  virtual element method}, ESAIM Math. Model. Numer. Anal., 50 (2016),
  pp.~879--904.

\bibitem{beirao2013basic}
{\sc L.~Beir\~{a}o~da Veiga, F.~Brezzi, A.~Cangiani, G.~Manzini, L.~D. Marini,
  and A.~Russo}, {\em Basic principles of virtual element methods}, Math.
  Models Methods Appl. Sci., 23 (2013), pp.~199--214.

\bibitem{da2013virtual}
{\sc L.~Beir\~{a}o~da Veiga, F.~Brezzi, and L.~D. Marini}, {\em Virtual
  elements for linear elasticity problems}, SIAM J. Numer. Anal., 51 (2013),
  pp.~794--812.

\bibitem{beirao2014hitchhiker}
{\sc L.~Beir\~{a}o~da Veiga, F.~Brezzi, L.~D. Marini, and A.~Russo}, {\em The
  hitchhiker's guide to the virtual element method}, Math. Models Methods Appl.
  Sci., 24 (2014), pp.~1541--1573.

\bibitem{beirao2016virtual}
\leavevmode\vrule height 2pt depth -1.6pt width 23pt, {\em Virtual element
  method for general second-order elliptic problems on polygonal meshes}, Math.
  Models Methods Appl. Sci., 26 (2016), pp.~729--750.

\bibitem{da2017divergence}
{\sc L.~Beir\~{a}o~da Veiga, C.~Lovadina, and G.~Vacca}, {\em Divergence free
  virtual elements for the {S}tokes problem on polygonal meshes}, ESAIM Math.
  Model. Numer. Anal., 51 (2017), pp.~509--535.

\bibitem{brenner1990nonconforming}
{\sc S.~C. Brenner}, {\em A nonconforming multigrid method for the stationary
  {S}tokes equations}, Math. Comp., 55 (1990), pp.~411--437.

\bibitem{brezzi2014basic}
{\sc F.~Brezzi, R.~S. Falk, and L.~D. Marini}, {\em Basic principles of mixed
  virtual element methods}, ESAIM Math. Model. Numer. Anal., 48 (2014),
  pp.~1227--1240.

\bibitem{cangiani2016nonconforming}
{\sc A.~Cangiani, V.~Gyrya, and G.~Manzini}, {\em The nonconforming virtual
  element method for the {S}tokes equations}, SIAM J. Numer. Anal., 54 (2016),
  pp.~3411--3435.

\bibitem{cangiani2017conforming}
{\sc A.~Cangiani, G.~Manzini, and O.~J. Sutton}, {\em Conforming and
  nonconforming virtual element methods for elliptic problems}, IMA J. Numer.
  Anal., 37 (2017), pp.~1317--1354.

\bibitem{crouzeix1973conforming}
{\sc M.~Crouzeix and P.-A. Raviart}, {\em Conforming and nonconforming finite
  element methods for solving the stationary {S}tokes equations. {I}}, Rev.
  Fran\c{c}aise Automat. Informat. Recherche Op\'{e}rationnelle S\'{e}r. Rouge,
  7 (1973), pp.~33--75.

\bibitem{liu2017nonconforming}
{\sc X.~Liu, J.~Li, and Z.~Chen}, {\em A nonconforming virtual element method
  for the {S}tokes problem on general meshes}, Comput. Methods Appl. Mech.
  Engrg., 320 (2017), pp.~694--711.

\bibitem{talischi2012polymesher}
{\sc C.~Talischi, G.~H. Paulino, A.~Pereira, and I.~F.~M. Menezes}, {\em {\tt
  {P}oly{M}esher}: a general-purpose mesh generator for polygonal elements
  written in {M}atlab}, Struct. Multidiscip. Optim., 45 (2012), pp.~309--328.

\bibitem{thomassetimplementation}
{\sc F.~Thomasset}, {\em Implementation of finite element methods for
  {N}avier-{S}tokes equations}, Springer Series in Computational Physics,
  Springer-Verlag, New York-Berlin, 1981.

\bibitem{zhang2019nonconforming}
{\sc B.~Zhang, J.~Zhao, Y.~Yang, and S.~Chen}, {\em The nonconforming virtual
  element method for elasticity problems}, J. Comput. Phys., 378 (2019),
  pp.~394--410.

\bibitem{zhao2019divergence}
{\sc J.~Zhao, B.~Zhang, S.~Mao, and S.~Chen}, {\em The divergence-free
  nonconforming virtual element for the {S}tokes problem}, SIAM J. Numer.
  Anal., 57 (2019), pp.~2730--2759.

\end{thebibliography}

\end{document}